\documentclass[11pt]{article} 
\usepackage[T1]{fontenc}
\usepackage[english]{babel}

\usepackage[utf8]{inputenc} 

\usepackage{geometry} 
\usepackage{graphicx} 

\usepackage{booktabs} 
\usepackage{array} 
\usepackage{paralist} 
\usepackage{verbatim} 
\usepackage{subfig} 

\usepackage{csquotes}
\usepackage{graphicx} 
\usepackage{amssymb}
\usepackage{amsmath, amsfonts, bbm, dsfont}
\usepackage[all]{xy}
\usepackage{MnSymbol}

\usepackage{fancyhdr} 
\usepackage{sectsty}
\allsectionsfont{\sffamily\mdseries\upshape} 

\usepackage[nottoc,notlof,notlot]{tocbibind} 
\usepackage[titles,subfigure]{tocloft} 

\usepackage[
backend=biber,
style=alphabetic,
giveninits=true,
doi=false,
url=false,
isbn=false
]{biblatex}
\renewbibmacro{in:}{}

\usepackage{amsthm}
\usepackage{tikz-cd}

\tikzcdset{row sep/normal= 0.8 cm,
	column sep/normal= 0.8 cm}

\DeclareMathOperator{\flex}{\textbf{\textasciicircum}}

\newtheorem{theo}{Theorem}[section]
\newtheorem{prop}[theo]{Proposition}
\newtheorem{lemma}[theo]{Lemma}
\newtheorem{coro}[theo]{Corollary}

\newtheorem{claim}[theo]{Claim}

\theoremstyle{definition}
\newtheorem{fact}[theo]{Fact}
\newtheorem{defi}[theo]{Definition}
\newtheorem{rem}[theo]{Remark}
\newtheorem{ex}[theo]{Example}
\newtheorem{contre-ex}[theo]{Contre-exemple}

\usepackage{hyperref}
\hypersetup{
	colorlinks=true,
	linkcolor=black,
	filecolor=black,      
	urlcolor=black,
	citecolor=black,
}

\usepackage{tikz}
\usetikzlibrary{intersections}

\title{On groups and fields interpretable in $\mathrm{NTP_2}$ fields}

\author{Paul $\mathrm{Wang}^1$}

\date{}

\begin{document}
\setcounter{tocdepth}{2}
\maketitle

  ABSTRACT.  \small This paper aims at developing model-theoretic tools to study interpretable fields and definably amenable groups, mainly in $\mathrm{NIP}$ or $\mathrm{NTP_2}$ settings. An abstract theorem constructing definable group homomorphisms from generic data is proved. It relies heavily on a stabilizer theorem of Montenegro, Onshuus and Simon. The main application is a structure theorem for definably amenable groups that are interpretable in algebraically bounded perfect $\mathrm{NTP_2}$ fields with bounded Galois group (under some mild assumption on the imaginaries involved), or in algebraically bounded theories of (differential) NIP fields. These imply a classification of the fields interpretable in differentially closed valued fields, and structure theorems for fields interpretable in henselian valued fields of characteristic $0$, or in NIP algebraically bounded differential fields.
   \footnotetext[1]{partially funded by ANR GeoMod  (AAPG2019, ANR-DFG)}
 \footnotetext[2]{\copyright \,  2024. This manuscript version is made available under the CC-BY-NC-ND 4.0 license http://creativecommons.org/licenses/by-nc-nd/4.0/}
 
\normalsize

\section{Introduction}

The study of algebraic structures, say, groups, that are definable in given first-order theories is a part of what is called \emph{geometric} model theory. In principle, the presence of definable groups or definable fields sheds light upon the model theory of the structure under study. One instance of this phenomenon is Hrushovski's proof of the fact that unidimensional theories are superstable (see \cite{Hru-UniDim} or \cite[Theorem 5.22]{poizat-stable-groups}), which goes through the study of definable groups, even if they are not mentioned in the statement.

Recall that, given a first-order structure $S$, one can define the \emph{category of definable sets} in $S$, where the objects are the definable sets (with parameters), and the morphisms are the definable maps, i.e. those whose graphs are definable sets. Adjoining all definable quotients, i.e. imaginaries, to the category of definable sets, one constructs the category of \emph{interpretable} sets. Definable groups, resp. interpretable groups, can be viewed as group objects in categories of definable sets, resp. interpretable sets. Studying intepretable groups can then be rephrased as trying to describe categories of interpretable groups. One particularly interesting case is that where $S$ contains a (possibly enriched) field structure $(K, +, \times)$. Then, the category of algebraic varieties over $K$ is a subcategory of the category of definable sets in $S$. Hence, algebraic groups over $K$ are also definable groups in $S$. In that context, one may try to find ways of \emph{comparing the category of definable groups}, or interpretable groups, \emph{to the subcategory of algebraic groups}.

Classical results in this area include the fact, commonly attributed to Weil-Hrushovski (see for instance \cite[Theorem 4.13]{poizat-stable-groups}), that any group definable in an algebraically closed field is definably isomorphic to an algebraic group, and the consequence, due to Poizat, that any infinite field definable in an algebraically closed field is definably isomorphic to the ambient field (see Theorem 4.15 in \cite{poizat-stable-groups}).

The case of algebraically closed fields is, from several points of view, the simplest one. The following results deal with the more subtle - and also more interesting to number theorists - cases of real, p-adic or pseudo-finite fields.

\begin{theo}[see Theorem A in \cite{HruPil-GpPFF}]

Let $F$ be either $\mathbb{R}$ or $\mathbb{Q}_p$. Then, any Nash group over $F$ is locally (i.e. in neighbourhoods of the identity) Nash isomorphic to the set of $F$-rational points of an algebraic group defined over $F$.

\end{theo}

\begin{theo}[see Theorem C in \cite{HruPil-GpPFF}]

Let $G$ be a group definable in a pseudo-finite field $F$. Then, $G$ is virtually isogeneous with the set of $F$-rational points of an algebraic group defined over $F$.

\end{theo}

The abstract model-theoretic result underlying these theorems is expressed in terms of ``geometric substructures of strongly minimal sets'' (see Definition 2.6 in \cite{HruPil-GpPFF}). The proof relies on a suitable application of the guiding principles of the famous \emph{group configuration theorem} (see for instance Chapter 5 of \cite{pillay}), namely, the reconstruction of a definable group from generic data.

Recently, several similar results have been obtained. For instance, in the context of valued fields of dp-rank $1$ (a tameness property which is stronger than NIP), under some technical assumptions, the classification of interpretable fields has been carried out in \cite[Theorem 7.1]{fields_interp_in_various_val_fields}. Other results include the work done in \cite{ppp_def_groups_codf}, which proves that groups definable in various differential fields can be definably embedded into algebraic groups, and \cite[Theorem 1]{semisimple_groups_interp_in_various_val_fields}, dealing with definably semisimple groups interpretable in several valued fields.

In this paper, we try to develop these ideas, following similar principles, but using distinct model-theoretic tools. While several of the cases we study here are also dealt with in the works cited above, our methods differ, and many new examples are covered (for instance, see Example \ref{ex_classes_of_examples}). Also, our approach can deal with imaginaries, which is not done in \cite{ppp_def_groups_codf}, but is in \cite{fields_interp_in_various_val_fields} and \cite{semisimple_groups_interp_in_various_val_fields}, although in a different way. We do not require dp-minimality, and some results even hold for $\mathrm{NTP}_2$ fields, which is a larger class than NIP (see Definition \ref{defi_NTP2}). 


%
%

Section 2 introduces the main model-theoretic tools we will use to study definable groups. First, we recall a few facts about ideals, f-genericity and definably amenable groups in $\mathrm{NTP_2}$. We then state a crucial technical result, the stabilizer Theorem \ref{theo_stabilizers}, which is \cite[Theorem 2.15]{Stabilizers-NTP2}. We also state, and prove, our variant (Theorem \ref{theo_relative_stable_group_config}) of the usual stable group configuration theorem. The ideas there are similar to \cite{HruPil-GpPFF}. Combining the two, we get the following abstract result :
 \begin{theo}[Theorem \ref{theo_group_with_kernel_blind_to_S}]\label{theo_1}
 
  Let $T_0, T_1$ be theories, in languages $\mathcal{L}_0 \subseteq \mathcal{L}_1$ respectively. Assume that $T_0$ is superstable, has quantifier elimination and elimination of imaginaries in $\mathcal{L}_0$, in a collection of sorts $\mathcal{S}$. Assume that $T_1$ is $\mathrm{NIP}$. Let $N_0 \models T_0$, $N_1 \models T_1$, and $\iota : N_1 \rightarrow N_0$ be an $\mathcal{L}_0$-embedding. Assume that $dcl_0(\iota(N_1)) = \iota(N_1)$ and $acl_0(A) \cap \iota(N_1) = acl_1(A)$ for all $A \subseteq \iota(N_1)$.
  
   Let $V$ be a product of sorts of $\mathcal{S}$. Let $G$ be a definably amenable group, definable in $N_1$, with an $N_1$-definable map $\pi: G \rightarrow V$, whose fibers are blind\footnotemark[1] to $\mathcal{S}$. Then, there exists a quantifier-free-$\mathcal{L}_0(N_1)$-definable group $H$ in the sorts of $\mathcal{S}$, and an $\mathcal{L}_1(N_1)$-definable group homomorphism $G^{00}_{N_1}(N_1) \rightarrow H(N_1)$ whose kernel is blind to the sorts in $\mathcal{S}$.
 
   \footnotetext[1]{A definable set $X$ is blind to a definable set $Y$ if all definable finite correspondences $X \rightarrow Y$ have finite image.}

 \end{theo}

One of the main points is that this theorem also tackles interpretable sets, and decomposes a group between imaginary and real sorts. Also, its general framework captures many examples, and is not actually restricted to the case of (possibly enriched) fields, even though we shall only use it for that purpose here.

In Section 3, we apply these tools to the case of theories of possibly enriched fields. There, we may replace the $\mathrm{NIP}$ requirement with $\mathrm{NTP_2}$, at the cost of additional algebraic hypotheses. The following theorem was inspired by \cite[Proposition 6.1]{HruRK-MetaGp} and \cite[Theorem 2.19]{Stabilizers-NTP2}, and follows from Theorem \ref{theo_1} above: 

\begin{theo}[Theorem \ref{theo_general_morphism_with_imaginary_kernel}]

Let $T$ be an algebraically bounded (see Definition \ref{defi_alg_bounded}) theory of possibly enriched, possibly many-sorted, perfect fields. Assume that $T$ is $\mathrm{NIP}$, or that it is $\mathrm{NTP_2}$ and its models have bounded Galois groups. Let $G$ be a definably amenable group, admitting a definable map to some $K^n$, whose fibers are blind to $K$, in some sufficiently saturated model $M$ of $T$. Then, there exist an algebraic group $H$ over $M$ and an $M$-definable homomorphism $G^{00}_M \rightarrow H$, whose kernel is purely imaginary.
\end{theo}

Note that, in this theorem, the hypothesis on $G$ is mild: in most examples, the imaginary sorts we consider (other than the field sort itself) are blind to the field sort.
Applying this theorem to groups of affine transformations, we then prove the following
\begin{theo}[Corollary \ref{coro_fields_are_either_pur_im_or_algebraic}]
Let $T$ be an algebraically bounded theory of possibly enriched, possibly many-sorted, perfect fields. Assume that $T$ is $\mathrm{NIP}$, or that it is $\mathrm{NTP_2}$ and its models have bounded Galois groups. Let $F$ be an infinite definable field, admitting a definable map to some $K^n$, whose fibers are blind to $K$, in some sufficiently saturated model $M$ of $T$. Then either $F$ is purely imaginary, or $F$ admits a definable embedding into some finite extension of $K(M)$.

\end{theo}

One key ingredient is a general result, namely Proposition \ref{prop_fields_whose_affine_group_embeds_in_an_algebraic_group_are_known}, which states that, if $F$ and $K$ are definable fields, such that the group of affine transformations of $F$ admits a definable embedding into an algebraic group over $K$, then $F$ can be embedded definably into some finite extension of $K$. This result only uses so-called ``model-theoretic algebra'', and does not require any tameness assumption on the ambient theory.

In Section 4, as an illustration of the potential of these abstract tools, we study several examples. The first is the class of NIP \emph{algebraically bounded differential fields}, i.e. those where the algebraic closure is tame (see Definition \ref{defi_alg_bounded_diff_fields}). The general tools of Section 2 yield the following
\begin{theo}[Theorem \ref{theo_groups_in_algebraically_bounded_differential_fields}]

Let $T$ be an NIP theory of algebraically bounded differential fields. Let $M$ be a sufficiently saturated model of $T$. Let $n < \omega$.
\begin{enumerate}
\item Let $G$ be an $M$-definable group, admitting an $M$-definable map $G \rightarrow K^n$ whose fibers are blind to $K$. Assume that $G$ is definably amenable. Then, there exists an $M$-definable group homomorphism $G^{00} \rightarrow H$ with kernel blind to $K$, where $H$ is an algebraic group over $K(M)$.

\item Let $F$ be an infinite $M$-definable field, admitting an $M$-definable map $F \rightarrow K^n$ whose fibers are blind to $K$.  Then, either $F$ is blind to $K$, or $F$ admits a definable embedding into some finite extension of $K(M)$.

\end{enumerate}

\end{theo}

Then, we study the case of $\mathrm{DCVF}$, the theory of algebraically closed valued fields of equicharacteristic $0$ with a generic derivation, i.e. the model completion of the theory of differential valued fields of equicharacteristic $0$, which belongs to that class (see Fact \ref{fact_DCVF}). In Theorem \ref{theo_fields_interpretable_in_DCVF}, we show that, in $\mathrm{DCVF}$, the only infinite interpretable fields are the residue field, the valued field, and the field of constants. As a second example, we prove the following theorem for Henselian valued fields of characteristic $0$:

\begin{theo}[Theorem \ref{theo_fields_hen_0}]
Let $T \supseteq Hen_0$ be a theory of finitely ramified (see Definition \ref{defi_fin_rami}) pure henselian valued fields of characteristic $0$ with perfect residue fields, in the three-sorted language $(K, k, \Gamma)$. Assume that $T$ is NIP, or that $T$ is $\mathrm{NTP_2}$ with $Gal(K^{alg}/K)$ bounded in all models of $T$. Also assume that $k$, with its induced structure, is algebraically bounded as a field. Let $M$ be a sufficiently saturated model of $T$. Let $m, n, l < \omega$. Let $F \subseteq K^n \times k^m \times \Gamma^l$ be an $M$-definable infinite field. Then, exactly one of the following holds:

\begin{enumerate}

\item The field $F$ admits a definable embedding into some finite extension of $K$.
 
 \item The field $F$ admits a definable embedding into some finite extension of $k$.
 
 \item The field $F$ is definably isomorphic to a field definable in $\Gamma$.

\end{enumerate}

\end{theo}

\medskip \medskip

In this paper, if $T$ is a fixed complete theory, $M$ is a model of $T$, and $A \subseteq M$, we say that $A$ is \emph{small} with respect to $M$ if $M$ is $|A|^+$-saturated and $|A|^+$-strongly homogeneous. We may write $A \subseteq^+ M$ to denote that smallness condition. Similarly, we may write $M \preceq^+ N$ if $M$ is a small elementary substructure of $N$.

\begin{section}{Ideals, f-generics, stabilizer theorems and group configurations}

	In this section, we describe in detail the model-theoretic tools we will apply in the study of definable groups and fields. We shall introduce a stabilizer theorem (Theorem \ref{theo_stabilizers}) and a group configuration theorem (Theorem \ref{theo_relative_stable_group_config}). Let us first recall the definitions of $\mathrm{NTP_2}$ and NIP.

	\begin{defi}\label{defi_NTP2}
		
		We say that a partitioned formula $\phi(x,y)$ has $\mathrm{TP_2}$, with respect to a theory $T$, if there are $(a_{l,j})_{l,j < \omega}$ in some $M \models T$ and $k < \omega$ such that:
		
		\begin{enumerate}
			
			\item The collection of formulas $\lbrace \phi(x, a_{l,j}) \, | \, j < \omega \rbrace$ is $k$-inconsistent for all $l < \omega$.
			
			\item For all $f : \omega \rightarrow \omega$, the collection of formulas $\lbrace \phi(x, a_{l, f(l)}) \, | \, l < \omega \rbrace$ is consistent.
			
		\end{enumerate} 
		
		A formula has $\mathrm{NTP_2}$ if it does not have $\mathrm{TP_2}$. The theory $T$ is $\mathrm{NTP_2}$ if no formula has $\mathrm{TP_2}$.

		Similarly, $\phi(x,y)$ has IP if there are $(a_i)_{i < \omega}$ and $(b_J)_{J \in P(\omega)}$ in some $M \models T$ such that for all $i$, $J$, we have $M \models \phi(a_i, b_J)$ if and only if $i \in J$. A formula is NIP if it does not have IP, and a theory is NIP if all formulas are NIP.
		
	\end{defi}

	At the level of theories, $\mathrm{TP_2}$ implies IP. In other words, all NIP theories are $\mathrm{NTP_2}$. See \cite[Proposition 5.31]{Sim-Book}

	Let us also recall some notions about connected components of type-definable groups.
	
	\begin{defi}

		Let $G$ be a type-definable group. Let $A$ be a set of parameters such that $G$ is $A$-type-definable. Then, we let $G^0_A$, resp. $G^{00}_A$ denote the intersection of the relatively $A$-definable, resp. $A$-type-definable, subgroups of $G$ of bounded index.

		We say that $G^0$ \emph{exists}, resp. $G^{00}$ exists, if the type-definable subgroup $G^0_A$, resp. $G^{00}_A$, does not depend on $A$ (as long as $G$ is $A$-type-definable). 
		The group $G$ is said to be \emph{connected} if $G=G^{00}_B$ for all $B$.
	\end{defi}

	\begin{fact}\label{fact_G^00}
		\begin{enumerate}
			\item In stable theories, $G^0$ and $G^{00}$ always exist, and they are equal. In $\omega$-stable theories, if $G$ is definable, then $G^0$ is definable and of finite index in $G$.
			
			\item In NIP theories, $G^{00}$ always exists. (see for instance \cite[Theorem 8.7]{Sim-Book}) 
			
			\item In general, for any $A$, the subgroup $G^{00}_A$ is a normal subgroup of $G$. See \cite[Lemma 4.1.11]{Wag-Simple}.
			
		\end{enumerate}
		
	\end{fact}

	\begin{subsection}{Ideals}

		In this subsection, we fix a complete theory $T$, which we assume to eliminate quantifiers and imaginaries, and a very saturated and very homogeneous model  $\mathcal{U}$. A subset $A \subseteq \mathcal{U}$ is \emph{small} if it is small with respect to $\mathcal{U}$, i.e. if the model $\mathcal{U}$ is $|A|^+$-saturated and $|A|^+$-strongly homogeneous.
		
		\begin{defi}
			
			\begin{enumerate}
				
				\item An \emph{ideal} $\mu$ is a collection of $\mathcal{U}$-definable sets, which is  closed under finite unions, and such that any $\mathcal{U}$-definable set contained in an element of $\mu$ is itself in $\mu$.  Let $A$ be a small set. The ideal $\mu$ is \emph{$A$-invariant} if, for all formulas $\phi(x,y)$ over $A$, and all elements $b, c \in \mathcal{U}$ such that $b \equiv_A c$, we have $\phi(x,b) \in \mu$ if and only if $\phi(x, c) \in \mu$.

				We say that a type-definable set (possibly a complete type) $\pi$ is in an ideal $\mu$ if some definable set $X$ containing $\pi$ is in $\mu$.
				
				\item If $I$ is a linearly ordered set, we say that a sequence of definable sets $(X_i)_{i \in I}$ is $A$-indiscernible if there exists a formula $\phi(x,y)$ over $A$, and an $A$-indiscernible sequence $(b_i)_{i \in I}$ such that, for all $i \in I$, the formula $\phi(x, b_i)$ defines $X_i$. This is equivalent to $A$-indiscernibility of the codes of the $X_i$.

				\item An $A$-invariant ideal $\mu$ has the \emph{S1 property} if, for any $A$-indiscernible sequence $(X_i)_{i < \omega}$ of definable sets, if $X_i \cap X_j$ is in $\mu$ for some/all $i \neq j$, then $X_i$ is in $\mu$ for some/all $i$.

				If $X$ is a definable set, we say that $\mu$ is $S1$ on $X$ if, the above property holds for $A$-indiscernible sequences $(X_i)_{i < \omega}$ of definable subsets of $X$, and if $X$ is not in $\mu$.

				If $\pi$ is a type-definable set, we will say that $\mu$ is $S1$ on $\pi$ if $\pi$ is not in $\mu$, and for some definable set $Y$ containing $\pi$, the ideal $\mu$ is $S1$ on $Y$.
				
			\end{enumerate}
			
		\end{defi}
		
		\begin{rem}
			Ideals correspond exactly to ring-theoretic ideals of the Boolean ring $Def_x(\mathcal{U})$ of  $\mathcal{U}$-definable sets (whose operations are symmetric difference and intersection).

			With this point of view, $A$-invariant ideals can also be thought of as families $(I_B)_{A \subseteq B \subseteq \mathcal{U}}$ such that, for all $B$, $I_B$ is an ideal of $Def_x(B)$, with the following conditions :
			
			\begin{enumerate}[(i)]
				\item If $A \subseteq B \subseteq C$, and $f : Def_x(B) \hookrightarrow Def_x(C)$ is the canonical embedding, then $I_B = f^{-1}(I_C)$.
				
				\item For all models $M \supseteq A$ contained in $\mathcal{U}$, the ideal $I_M \subseteq Def_x(M)$ is fixed setwise by the group $Aut(M/A)$.
			\end{enumerate}

			With this point of view, one could define $A$-invariant ideals without any reference to a saturated model $\mathcal{U}$. Note that the notion of $A$-invariant ideal can be used to generalize the notion of $A$-invariant complete type, applying negations. Unless otherwise mentioned, all ideals are proper.
			
		\end{rem}
		
		\begin{ex}
			Let $m$ be a finitely additive positive measure on definable sets, which is $A$-invariant. Then, the collection $\mu$ of measure-zero definable sets is an $A$-invariant ideal. Moreover, if $X$ is a definable set of finite positive measure, then $\mu$ is $S1$ on $X$.  
		\end{ex}

		\begin{rem}
			In \cite{Stabilizers-NTP2}, sets which are not in a given ideal $\mu$ are called $\mu$-\emph{wide}. This terminology seems slightly misleading : in the example of a measure, the ideal is that of measure-zero sets. Then, the wide definable sets would be those of positive measure, not those of full measure. Also, inspecting the definitions, a type-definable set could be wide and have measure zero: for instance, in RCF, the type-definable subgroup of infinitesimals has measure zero for the nonstandard Lebesgue measure, however any definable set containing it has positive measure.

			On the other hand, the sets on which an invariant ideal $\mu$ is $S1$ are called \emph{medium} sets in \cite{Stabilizers-NTP2}, and this choice seems reasonable.
		\end{rem}

		There is an interesting connection between forking and S1 ideals, which follows from the definitions:
		
		\begin{fact}[See Lemma 2.9 in \cite{Hru-StaApp}]\label{fact_S1_and_forking}
			Let $M \models T$, and $\mu$ an $M$-invariant ideal. Let $B \supseteq M$, and $p \in S(B)$ which is not in $\mu$, such that $\mu$ is S1 on $p|_M$. Then $p$ does not fork over $M$.
		\end{fact}

		The following notion is a generalization of what was called ``being purely imaginary'' in \cite{HruRK-MetaGp}.
		
		\begin{defi}
			A definable \emph{finite correspondence} $X \rightarrow Y$ is a definable subset $R \subseteq X \times Y$ such that, for all $x \in X$, the fiber $R_x \subseteq Y$ is finite, and non empty.
			
			We say that a definable set $X$ is \emph{blind} to a definable set $Y$ if all definable finite correspondences $X \rightarrow Y$ have finite image. If $\mathcal{Y}$ is a collection of definable sets, we say that $X$ is blind to $\mathcal{Y}$ if $X$ is blind to all the definable sets in $\mathcal{Y}$.

			We say that two definable sets are \emph{orthogonal} if they are blind to each other.
			
		\end{defi}

		We note that the collection of definable sets which are blind to a given definable set $Y$ is an invariant ideal. In fact, we can say a bit more :

		\begin{defi}
			Let $A$ be small set of parameters, and $\lambda$ be an $A$-invariant ideal. We say that it admits \emph{non-forking descent} if, for any $C \supseteq B \supseteq A$, for any tuple $a$ (in the appropriate sorts), if $a \downfree_B C$ or  $C \downfree_B a$ and $tp(a/C)$ is in $\lambda$, then $tp(a/B)$ is in $\lambda$.
			
			We say that it is \emph{ind-definable}, or $\vee$-definable, if, for any definable family $(X_u)_{u \in U}$ of definable sets, the collection of $u \in U$ such that $X_u \in \nu$, is $\vee$-definable, i.e. corresponds to an open set in the relevant type space. Note that this is weaker than definability.
		\end{defi}

		\begin{rem}
			Let $A$ be some set of parameters, and let $\lambda_A$ denote the $A$-invariant ideal of definable sets which fork over $A$. If the ideal $\lambda_A$ admits nonforking descent, then, for all $C \supseteq B \supseteq A$, for all $q \in S(C)$, if $q$ does not fork over $B$ and $q|_B$ does not fork over $A$, then $q$ does not fork over $A$.
			
		\end{rem}

		\begin{prop}\label{prop_the_blinders_are_a_nice_ideal}
			Let $\mathcal{Y}$ be a collection of definable sets. Let $I$ be a definable set blind to $\mathcal{Y}$, and $V$ be a product of sets in $\mathcal{Y}$. Let $\nu$ be the collection of definable subsets of $V \times I$ which are blind to $\mathcal{Y}$. Then, $\nu$ is an invariant, $\vee$-definable ideal, which is closed under definable bijections and admits non-forking descent.
			
		\end{prop}
		
		\begin{proof} From the definitions, it follows that $\nu$ is an ideal, closed under definable bijections. Also, if $A$ is a set of parameters over which everything is defined, then $\nu$ is $A$-invariant. The only facts left to check are the $\vee$-definability and non-forking descent. 
			
			\begin{claim}
				Let $X \subseteq V \times I$ be definable. Then, $X$ is in $\nu$ if and only if its projection to $V$ is finite.
			\end{claim}
			\begin{proof}
				If $X$ is in $\nu$, then it is blind to $\mathcal{Y}$. Since $V$ is a product of sets in $\mathcal{Y}$, the projection $X \rightarrow V$ has finite image. Conversely, if $\pi(X)$ is finite, then $X$ is in definable bijection with a finite union of subsets of $I$. Since $I$ is blind to $\mathcal{Y}$, so is $X$.
			\end{proof}
			
			The claim implies that $\nu$ is $\vee$-definable, since finiteness of the projection to $V$ is $\vee$-definable. 
			
			Let us now prove nonforking invariance. Let $C \supseteq B \supseteq A$ be parameter sets, and $a = (a_V, a_I) \in V \times I$. Let us assume that $tp(a / C)$ is in $\nu$. Then, from the above, we have $a_V \in acl(C)$. Now, if either $a \downfree_B C$ or $C \downfree_B a$, then $a_V \in acl(B)$, i.e. $tp(a/B) \in \nu$. This shows non-forking descent, and concludes the proof. \end{proof}

		%

		Let us conclude this section by spelling out two lemmas on ideals that are implicit in \cite{Stabilizers-NTP2}, and that we will use later.

		\begin{lemma}\label{lemma_S1_and_finite_fibers}
			Let $\mu_1$ be an S1 and $M$-invariant ideal on some $M$-definable set $A$. Let $B$ be another $M$-definable set. Let $\mu$ be the ideal on $A \times B$ containing all the sets whose projections to $A$ are in $\mu_1$.
			
			Let $X \subseteq A \times B$ be such that the projection $\pi_A: X \mapsto A$ has finite fibers. Then, either $X$ is in $\mu$, or $\mu$ is S1 on $X$.

		\end{lemma}

		\begin{proof}
			Let us assume that $X$ is not in $\mu$, and show that $\mu$ is S1 on $X$. Let $(X_i)_{i < \omega}$ be an $M$-indiscernible sequence of definable subsets of $X$, such that $X_i \cap X_j \in \mu$, for all $i \neq j$. We need to show that some (equivalently, all) $X_i$ are in $\mu$. So, we consider the sets $\pi_A(X_i)$, and we want to show that they are in $\mu_1$.

			By compactness, there is a uniform bound $k < \omega$ on the size of fibers, in $X$, of elements of $\pi_A(X)$.
			
			\begin{claim} For all $i_1 < i_2 < \cdots < i_{k+1}$, we have the following inclusion:  $\pi_A(X_{i_1}) \cap \cdots \cap \pi_A(X_{i_{k+1}}) \subseteq \bigcup\limits_{1 \leq r < s \leq k+1} \pi_A(X_{i_r} \cap X_{i_s})$. 
			\end{claim}
			\begin{proof}
				We use the pigeonhole principle. Let $a \in \pi_A(X_{i_1}) \cap \cdots \cap \pi_A(X_{i_{k+1}})$. So, for $j=1, \cdots, k+1$, the element $a$ admits a preimage in $X_{i_j}$. However, the fiber over $a$ has size at most $k$, so one of these preimages belongs to $X_{i_j} \cap X_{i_{l}}$ for $l \neq j$. This proves the inclusion above.
				
			\end{proof}
			
			Recall that the $X_i \cap X_j$ are in $\mu$, i.e. $\pi_A(X_i \cap X_j) \in \mu_1$, for $i \neq j$. Thus, the set $\pi_A(X_{i_1}) \cap \cdots \cap \pi_A(X_{i_{k+1}})$ is in the ideal  $\mu_1$. We may assume that $k+1$ is a power of $2$. Then, since $\mu_1$ is S1 and $(\pi_A(X_i))_i$ is an indiscernible sequence of definable sets, it is straightforward to deduce by induction that $\pi_A(X_i)$ is in $\mu_1$ for all $i$, as required.
		\end{proof}

		\begin{lemma}\label{lemma_wide_over_acl}
			Let $A \subseteq B$ be small parameter sets, let $\mu$ be an $A$-invariant ideal, and let $a$ be an element such that $tp(a / B)$ is not in $\mu$. Then $tp(a / acl(B))$ is not in $\mu$.

		\end{lemma}
		
		\begin{proof}
			Let $X$ be an $acl(B)$-definable set such that $a \in X$. Let $(X_i)_{i <n}$ be the $B$-conjugates of $X$. Then, the set $\bigcup_{i <n} X_i$ is $B$-definable and contains $a$, so is not in $\mu$. Moreover, by $A$-invariance of $\mu$, either each $X_i$ is in $\mu$, or none of them are. Since $\mu$ is an ideal, it is the latter. \end{proof}

	\end{subsection}

	\begin{subsection}{Forking-genericity, definably amenable $\mathrm{NTP_2}$ groups and stabilizers}
		
		In this subsection, we shall work inside an $\mathrm{NTP_2}$ theory $T$. If $G$ is a definable group, and $g,h$ are elements of $G$, we let $gh$ denote the concatenation, or the union, of the tuples, and $g \cdot h$ the product.
		
		\begin{defi}
			A set of parameters $A$ is an extension base if no type over $A$ forks over $A$. 
		\end{defi}
		
		\begin{defi}
			
			Let $A$ be a set of parameters, $M$ a model containing $A$, $G$ an $A$-type-definable group, and $p \in S_G(M)$. 
			
			\begin{enumerate}
				\item We define a left action of $G(M)$ on $S_G(M)$ in the following way: for $g \in G(M)$, we define $g\cdot p$ as $tp(g \cdot a / M)$ for any $a$ realizing $p$. This does not depend on the choice of $a$. We also define a right action of $G(M)$, in the obvious way.

				\item If $M$ is $|A|^{+}$-saturated, we say that $p$ is strongly (left) f-generic over $A$ if, for all $g \in G(M)$, the type $g \cdot p$ does not fork over $A$.

				Similarly, $p$ is strongly bi-f-generic over $A$ if, for all $g,h \in G(M)$, the type $g \cdot p \cdot h$ does not fork over $A$.

			\end{enumerate}
			
		\end{defi}

		If $M, N$ are models of $T$, we write $M \preceq^+ N$ to say that $N$ is an $|M|^+$-saturated and $|M|^+$-strongly homogeneous elementary extension of $M$.

		\begin{fact}[See Lemma 3.11 in \cite{CS18}]
			Let $A \subset M \preceq N$, where $M$, $N$ are $|A|^{+}$-saturated models. Let $G$ be an $A$-definable group. Let $p \in S_G(M)$ be strongly (left, right, or bi) f-generic over $A$. Then, there exists $q \in S_G(N)$ which extends $p$ and is also strongly (left, right, or bi) f-generic over $A$. In particular, $q$ is a nonforking extension of $p$.
		\end{fact}

		For the following facts, see \cite[Corollary 5.2, Proposition 5.3, Proposition 5.9, Corollary 5.19]{poizat-stable-groups} and \cite[Definition 7.2.1, Corollary 7.2.4, Lemma 7.2.6]{Mar-IntroMT}, or \cite[Lemmas 3.11 and 3.12]{HruRK-MetaGp}.
		
		\begin{rem}\label{rem_generics_in_stable_groups}
			If $G$ is a type-definable group in a stable theory, the notion of f-genericity has some nice properties:

			\begin{enumerate}
				\item Types that are f-generic always exist, and they are characterized by boundedness of their orbit, and many other criteria. Usually, they are simply called ``generics''. Any left or right translate of a generic is generic, and generics are all left and right translates of each other. The connected component $G^0$ has a unique generic type, called the principal generic.
				
				\item Any element of $G$ is the product of two generics. Any element in $G^0$ is the product of two realizations of the principal generic.
				
				\item Genericity for types is a closed notion, thus it makes sense to say that a definable subset of $G$ is generic. In the $\omega$-stable case, a definable subset of $G$ is generic if and only if it has the same dimension as $G$. Also, over any base, there are only finitely many generics.

			\end{enumerate}

		\end{rem}

		\begin{defi}\label{defi_St_and_Stab} Let $M$ be a model of $T$, $A \subseteq^+ M$ a small set of parameters, and let $G$ be an $A$-definable group.
			Let $q,r \in S_G(M)$, and $\mu$ an $A$-invariant ideal, closed under left translations. If neither $q$ nor $r$ is in $\mu$, we define $St_{\mu}(q,r)$ to be the following $M$-invariant subset of $G$: $St_{\mu}(q,r) = \lbrace g \in G \, | \, g \cdot q \cap r \notin \mu \rbrace$. Here, $g \cdot q \cap r$ is a type-definable set, and by definition, it is in $\mu$ if \emph{some} definable set containing it is in $\mu$.

			We also define $St_{\mu}(q) = St_{\mu}(q,q)$, and $Stab_{\mu}(q)$ as the subgroup of $G$ generated by $St_{\mu}(q)$. Note that $St_{\mu}(q)$ is closed under taking inverses.

			Similarly, using the right action of $G$ on the types, we define ${St_{r, \mu}}(q,r)$, ${St_{r, \mu}}(q)$, and ${Stab_{r, \mu}}(q)$. If the context is clear, we may drop the subscript $\mu$.

		\end{defi}

		The following proposition was implicitly in Section 3 of \cite{Stabilizers-NTP2}. We give the details here.

		\begin{prop}\label{prop_mu_A_is_S1}
			Let $M$ be a model of $T$, and $A \subseteq^+ M$ be an extension base. Let $G$ be an $A$-definable group. Assume that $G$ admits a type $p \in S_G(M)$ which is strongly bi-f-generic over $A$. Let $\mu_A$ denote the ideal of definable sets which do not extend to global types bi-f-generic over $A$. Then, the ideal $\mu_A$ is S1.
		\end{prop}

		\begin{lemma}\label{lemma_phi(gxb)_forks}
			Let $A$ be a set of parameters, and let $G$ be an $A$-definable group. Let $\phi(x, b)$ be a formula over $Ab$, which forks over $A$. Let $g,h$ be elements of $G$ such that $tp(g,h /Ab)$ does not fork over $A$. Then, the formula $\phi(g\cdot x \cdot h, b)$ forks over $A$.
		\end{lemma}
		
		\begin{proof}
			Assume that $\phi(g\cdot x \cdot h,b)$ does not fork over $A$. Let $a$ be an element of $G$ such that $tp(a/ Ab gh)$ does not fork over $A$ and $\models \phi(g\cdot a\cdot h, b)$. Then, we have $a \downfree_{A gh} b$, so $gah \downfree_{A g h} b$. Since we assumed that $g h \downfree_A b$, we deduce by transitivity $gah \downfree_A b$, in particular $g \cdot a \cdot h \downfree_A b$. Since $g\cdot a\cdot h$ realizes the formula $\phi(x,b)$, the latter does not fork over $A$.
		\end{proof}

		\begin{prop}\label{prop_universal_witness_for_bigenericity}
			
			Let $A$ be an extension base, $B \supseteq A$, $\phi(x)$ a formula over $B$. Let $M \supseteq^+ B$ be a model. Let $G$ be an $A$-definable group, and let $q_1$, resp. $q_2$, be a type in $S_G(M)$ which is strongly left-generic over $A$, resp. strongly right-generic over $A$. Let $g \in M$ realize $q_1|_B$, and $h$ realize $q_2|_{Bg}$. Then, the formula $\phi(x)$ extends to a type $p \in S_G(M)$ which is strongly bi-f-generic over $A$ if and only if the formula $\phi(g\cdot x\cdot h) =: g^{-1} \cdot \phi \cdot h^{-1}$, does not fork over $A$.
			
		\end{prop}
		
		\begin{proof}
			One implication is clear : if $\phi$ extends to a type $p \in S_G(M)$ strongly bi-f-generic over $A$, then no bi-translate of $\phi$ forks over $A$. Let us prove the other implication, using contraposition. So, assume that $\phi$ does not extend to a type $p \in S_G(M)$ strongly bi-f-generic over $A$. Then, the partial type $\lbrace \phi(x) \rbrace \cup \lbrace \neg \psi(g_1 \cdot x \cdot h_1, m) \, | \,g_1, h_1 \in G(M), \psi(x,m)$ is over $M$ and forks over $A \rbrace$ is inconsistent. Hence, by compactness, there is an integer $N$, formulas $\psi_i(x, m)$ that are over $M$ and which fork over $A$, and elements $g_i, h_i \in G(M)$, for $i< N$, such that  $\phi(x) \models \bigvee_i \psi_i(g_i \cdot x \cdot h_i, m)$. So, we have $\phi(g\cdot x\cdot h) \models  \bigvee_i \psi_i(g_i \cdot g  \cdot x \cdot  h \cdot  h_i, m)$. Note that, since the statement depends only on the type of $(g,h)$ over $B$, we may assume that $g$ realizes $q_1|_{B m (g_i  h_i)_{i<N}}$, and $h$ realizes $q_2|_{B m g (g_i  h_i)_{i<N}}$.

			Now, we wish to show that, for all $i$, the formula $\psi_i(g_i  \cdot g \cdot  x \cdot  h \cdot  h_i, m)$ forks over $A$.
			
			\begin{claim}
				We have, for all $i < N$, the following independence relation: 
				
				$g_i \cdot g,  h \cdot h_i \downfree_A Bm$.
			\end{claim}
			
			\begin{proof}
				First, by assumption of left-genericity of $q_1$, we know that $tp(g_i \cdot g / Bm (g_j,  h_j)_{j<N}) $ does not fork over $A$, for all $i<N$. In particular, we have $g_i \cdot g \downfree_A Bm$. Similarly, $tp(h \cdot h_i / g Bm (g_j,  h_j)_{j<N})$ does not fork over $A$, which implies $h \cdot h_i \downfree_A g_i \cdot g Bm$, then $h \cdot h_i \downfree_{A g_i \cdot g} Bm$.  Thus, by left transitivity for nonforking, $tp(g_i \cdot g,  h \cdot h_i / Bm)$ does not fork over $A$.
			\end{proof}
			
			Hence, we may apply Lemma \ref{lemma_phi(gxb)_forks}, to show that, for all $i < N$, the formula $\psi_i(g_i  \cdot g  \cdot x  \cdot h  \cdot h_i, m)$ forks over $A$. Recall that we have $\phi(g\cdot x\cdot h) \models  \bigvee_i \psi_i(g_i \cdot g  \cdot x \cdot  h \cdot  h_i, m)$. Therefore, the formula $g^{-1} \cdot \phi(x) \cdot h^{-1}$, which is the same as $\phi(g\cdot x \cdot h)$, forks over $A$. This concludes the proof. 
		\end{proof}

		\begin{proof}[Proof of Proposition \ref{prop_mu_A_is_S1}]
			Let $(a_i)_{i < \kappa}$ be an $A$-indiscernible sequence, which we may assume to lie inside $M$, and $\phi(x,y)$ an $A$-formula, such that $\phi(x, a_i) \notin \mu_A$ for all $i < \kappa$. Let $q$ be type in $S_G(M)$, strongly bi-f-generic over $A$. Let $g$ realize $q|_{A (a_i)_{i < \kappa}}$, and $h$ realize $q|_{A g (a_i)_{i < \kappa}}$. Then, enlarging $\kappa$, extending the sequence $(a_i)$ if necessary, by Erdös-Rado, there exists a subsequence $(a_j)_{j< \omega}$ which is indiscernible over $A g  h$. Then, by Proposition \ref{prop_universal_witness_for_bigenericity}, the formula $g^{-1 }\cdot \phi(x, a_j) \cdot h^{-1}$ does not fork over $A$, for all $j < \omega$.

			Since the forking ideal over an extension base is S1 in $\mathrm{NTP_2}$ theories (see \cite{BYChe}, Corollary 2.10), this implies that the formula $g^{-1 }\cdot (\phi(x, a_j) \wedge \phi(x, a_k)) \cdot h^{-1}$ does not fork over $A$, for all $j < k < \omega$. By Proposition \ref{prop_universal_witness_for_bigenericity} again, we deduce that $\phi(x, a_j) \wedge \phi(x, a_k)$ is not in $\mu_A$, for all $j<k<\omega$.
		\end{proof}

		\begin{defi}\label{defi_definably_amenable}
			
			A definable group $G$ is \emph{definably amenable} if, for some model $M$ of $T$, there is a left-invariant finitely additive probability measure (also known as Keisler measure) on the $M$-definable subsets of $G$. (For more on definably amenable groups, see \cite[Chapter 8]{Sim-Book}) 
			
		\end{defi}

		\begin{rem}\label{rem_def_amenable_groups}
			\begin{enumerate}
				\item It happens that, in the definition above, the existence of a left-invariant Keisler measure does not depend on the choice of the model $M$. 
				
				\item Any definable group which is amenable as a discrete group is a fortiori definably amenable. In particular, any solvable and definable group is definably amenable. On the other hand, any stable group is definably amenable.
				
			\end{enumerate}
		\end{rem}

		In this paper, we will mainly be interested in definably amenable $\mathrm{NTP_2}$ groups due to the following result :
		
		\begin{prop}\label{prop_exists_generics}
			Let $G$ be a definably amenable $\mathrm{NTP_2}$ group. Let $A$ be an extension base, and $M_1 \supseteq A$ an $|A|^{+}$-saturated model. Then, there exists $p \in S_G(M_1)$ which is strongly bi-f-generic over $A$.
		\end{prop}
		
		\begin{proof}
			By Proposition 3.20 in \cite{Stabilizers-NTP2}, there exists a model $M$ and a type $q \in S(\mathcal{U})$ which is strongly f-generic over $M$. Then, by Lemmas 3.4 and 3.5 in \cite{Stabilizers-NTP2}, there exists a type $p \in S_G(M_1)$ which is strongly bi-f-generic over $A$.
		\end{proof}

		All the results above on f-generics in $\mathrm{NTP_2}$ can be combined with the following theorem, to build definable group homomorphisms, in $\mathrm{NTP_2}$ settings.

		\begin{theo}[\cite{Stabilizers-NTP2}, Theorem 2.15]\label{theo_stabilizers}
			Let $G$ be a group definable in an arbitrary structure $M$, which is sufficiently saturated. Let $\lambda, \mu$ be $M$-invariant ideals on $G$, closed under left and right multiplication. Assume that, for all definable sets $X \in \lambda$, either $X$ is in $\mu$ or $\mu$ is $S1$ on $X$.

			Let $p \in S_G(M)$ be a type which is in $\lambda$, but not in $\mu$. Assume the following :
			
			\textbf{(A)} For all $q, r \in S_G(M)$, for all $c \models q$, $d \models r$ such that $d \downfree_M c$, if $tp(c\cdot d / M)$ or $tp(d \cdot c / M)$ is in $\lambda$, then $q$ is in $\lambda$.

			\textbf{(B)} For all $a$, $b$ realizing $p$ such that $b \downfree_M a$, the type $tp(a^{-1} b / M)$ is in $\lambda$.

			\textbf{(F)} There exist $a$, $b$ realizing $p$ such that $b \downfree_M a$ and $a \downfree_M b$.

			\smallskip
			
			\noindent Then, $Stab_{\mu}(p) = St_{\mu}(p)^2 = (pp^{-1})^2$ is a connected type-definable group, which is in $\lambda$ and not in $\mu$. Also, $Stab_{\mu}(p) \setminus St_{\mu}(p)$ is contained in a union of $M$-definable sets that are in $\mu$.

		\end{theo}
		
		From this theorem, we can deduce the following general comparison result :

		\begin{theo}\label{theo_how_to_find_group_homomorphisms}
			Let $T$ be an $\mathrm{NTP_2}$ theory. Let $M$ be a sufficiently saturated model, and $G, H$ be $M$-definable groups. Assume that $G$ has a strongly bi-f-generic type $p \in S_G(M)$. Let $\nu$ be a $\vee$-definable $M$-invariant ideal on $G$, invariant under left and right translations, and admitting non-forking descent.

			Assume that there are elements $a \models p|_M$, $b \models p|_{Ma}$, $\alpha, \beta \in H$ with the following properties :

			\begin{enumerate}
				\item We have $\alpha \in acl(M a)$, $\beta \in acl(M b)$ and $\beta \cdot \alpha \in acl(M b \cdot a)$.
				
				\item The type $tp(a / M \alpha)$ is in $\nu$.
				
			\end{enumerate}

			Then, there exists a finite subgroup $H_0 \leq H(M)$ and an $M$-definable group homomorphism $G^{00}_M \rightarrow H_1/H_0$, where $H_1 \rhd H_0$ is the centralizer of $H_0$, whose kernel is in the ideal $\nu$.

			More precisely, let $q_1$ denote $tp(a, \alpha/ M)$. Let $\mu_G$ denote the ideal on $G$ containing the definable sets which do not extend to a strongly bi-f-generic type over a sufficiently saturated model. Let $\mu$ the pullback of $\mu_G$ to $G \times H$. Then the graph of the homomorphism is given by the image of the subgroup $Stab_{\mu}(q_1) = (q_1 q_1^{-1})^2 \leq G \times H_1$ in $G \times (H_1 / H_0)$.
			
		\end{theo}

		\begin{proof}
			We wish to apply Theorem \ref{theo_stabilizers} above, using ideas similar to those in the proof of \cite[Theorem 2.19]{Stabilizers-NTP2}. First, the ideal $\mu_G$ is $M$-invariant, invariant under left and right translations, and S1 by Proposition \ref{prop_mu_A_is_S1}. Then, by definition $\mu$ is the ideal on $G \times H$ such that $X \in \mu$ if and only if $\pi_G(X) \in \mu_G$. Let  $\lambda$ be the ideal on $G \times H$ of the sets $X$ such that $\pi_G : X \rightarrow G$ has finite fibers and $\pi_H : X \rightarrow H$ has fibers in $\nu$. Finally, let $q_1 = tp(a,  \alpha / M)$, $q_2 = tp(b,  \beta / M)$ and $q_3 = tp(b \cdot a,  \beta \cdot \alpha / M)$. For the rest of the proof, the stabilizers we consider are taken with respect to the ideal $\mu$ only.

			We wish to apply the stabilizer theorem, i.e. Theorem \ref{theo_stabilizers}, to the data ($q_1$, $\mu$, $\lambda$) in the group $G \times H$. To do this, we have to check the hypotheses of said theorem. First, the fact that, for any set $X \in \lambda$, either $X$ is in $\mu$ or $\mu$ is S1 on $X$, follows from Lemma \ref{lemma_S1_and_finite_fibers}.
			
			By construction, the type $q_1=tp(a, \alpha / M)$ is not in $\mu$, because $tp(a / M) = p$ is not in $\mu_G$.

			\begin{claim}
				The type $q_1$ is in $\lambda$.
				
			\end{claim}
			\begin{proof}
				The element $\alpha$ is in $acl(Ma)$, so the ``finite fibers'' part holds. Let us check that $(a, \alpha)$ belongs to an $M$-definable set whose projection to $H$ has fibers  in $\nu$. We know that $tp(a /M \alpha) \in \nu$. Since $\nu$ is $\vee$-definable, there are $M$-definable sets $X \subseteq G \times H$ and $Y \subseteq H$, containing $(a, \alpha)$ and $\alpha$ respectively, such that, for all $y \in Y$, the fiber $X_y$ is in $\nu$. In other words, the fibers of the set $X \cap G \times Y$ over tuples in $H$ are in $\nu$, as desired. \end{proof}
			
			Note that condition (F) is automatically satisfied, due to $\mathrm{NTP}_2$ : the model $M$ is an extension base, so, by Proposition 3.7 of \cite{CheKap}, condition (F) holds for $q_1 \in S(M)$. Let us now check condition (A) : let $q,r \in S_{G \times H}(M)$, let $c \models q$, $d \models r$, with $d \downfree_M c$ and such that $tp(c\cdot d / M)$ or $tp(d \cdot c / M)$ is in $\lambda$. Let us show that $q$ is in $\lambda$. We shall deal with the case where  $tp(c\cdot d / M)$ is in $\lambda$, the other one being similar. Let us write $c = (c_1, c_2) \in G \times H$ and $d = (d_1, d_2) \in G \times H$. From the definition of $\lambda$, there are two things to prove : first, that $c_2$ is algebraic over $M c_1$, and then, that $tp(c_1 /M c_2)$ is in $\nu$.

			Let us check the first point : we know that $c_2 \cdot d_2 \in acl(M c_1 \cdot d_1)$ and $d_1 d_2 \downfree_M c_1 c_2$. So $c_2 \in acl(M c_1 \cdot d_1  d_2)$ and $d_1 d_2 \downfree_{M c_1} c_2$. Thus, $c_2 \in acl(M c_1)$, as required.

			Now, let us prove the second point. We know that $tp(c_1 \cdot d_1 / M c_2 \cdot d_2) \in \nu$. A fortiori, the type $tp(c_1 \cdot d_1 / M d_1,  c_2 \cdot d_2) $ is in $\nu$. So, by translation invariance, we have $tp(c_1 / M  d_1, c_2 \cdot d_2) \in \nu$. Recall that $d_1 d_2 \downfree_M c_1 c_2$, which implies $d_1 d_2 \downfree_{M c_2} c_1$. The conclusion then follows from the hypothesis that $\nu$ admits non-forking descent.

			%

			For condition (B), we use the same ideas as in \cite[Theorem 2.19]{Stabilizers-NTP2} : 
			
			\begin{claim} The $M$-invariant sets $St_r(q_i)$ are in $\lambda$, for $i=1,2,3$. 
			\end{claim}
			
			\begin{proof}
				Let $i$ be $1$, $2$ or $3$. Let $r \in S(M)$ be in the $M$-invariant set $St_r(q_i)$. So, by definition of $St_r(q_i)$, for some/any realization $c$ of $r$, the type-definable set $q_i \cdot c \cap q_i$ is not in $\mu$. Let us pick some such $c$. Then, by compactness, there exists a complete type $r^{\prime} \in S(Mc)$ belonging to the type-definable set $q_i \cdot c \cap q_i$, and which is not in $\mu$. Then, let $d \models r^{\prime}$. In particular, we have $d \models q_i$, the type $tp(d / Mc)$ is not in $\mu$, and $d \cdot c \models q_i$. So, as $\mu$ is S1 on the type $q_i = tp(d / M)$, we know by Fact \ref{fact_S1_and_forking} that $tp(d / M c)$ does not fork over $M$, i.e. $d \downfree_M c$. Also, the type $tp(d \cdot c / M) = q_i$ is in $\lambda$. So, by condition (A), the type $r$ is in $\lambda$, as required. 
			\end{proof}

			\begin{claim} The type $q_1$ is in $St_r(q_2, q_3)$.
			\end{claim}
			
			\begin{proof}
				We know that $(b, \beta) \models q_2$ and $(b, \beta) \cdot (a, \alpha) = (b \cdot a, \beta \cdot \alpha)$, where $(a, \alpha) \models q_1$ and $ (b \cdot a, \beta \cdot \alpha) \models q_3$. In other words, the element $ (b \cdot a, \beta \cdot \alpha)$ is in $q_3 \cap q_2 \cdot (a, \alpha)$. It remains to check that $tp (b \cdot a, \beta \cdot \alpha / M a \alpha)$ is not in $\mu$. Recall that $b \models p|_{Ma}$, so $tp(b / M a)$ is not in $\mu_G$. So, by right-translation invariance, the type $tp(b \cdot a / M a)$ is not in $\mu_G$. Then, by Lemma \ref{lemma_wide_over_acl} applied to the $M$-invariant ideal $\mu_G$, the type $tp(b \cdot a / acl(M a))$ is not in $\mu_G$, so in particular $tp(b \cdot a / M a \alpha)$ is not in $\mu_G$. Then, by definition of the ideal $\mu$ on $G \times H$, the type $tp (b \cdot a, \beta \cdot \alpha / M a \alpha)$ is not in $\mu$, as required.
			\end{proof}

			\begin{claim}
				For any $c$, $d$ realizing $q_1$ such that $d \downfree_M c$, we have $c \cdot d^{-1} \in St_r(q_2)$. In particular, $tp(c \cdot d^{-1} / M)$ is then in $\lambda$.

			\end{claim}
			
			\begin{proof}
				It is Lemma 2.9 in \cite{Stabilizers-NTP2}. \end{proof}

			So, condition (B) holds. Thus, we can finally apply the stabilizer theorem (Theorem \ref{theo_stabilizers}) to $\mu, \lambda, q_1$ in the group $G \times H$: the subgroup $Z = Stab(q_1) = (q_1 q_1^{-1})^2$ is $M$-type-definable, included in $\lambda$, connected, and is not in $\mu$. It remains to define a group morphism using $Z$.

			From the construction, the group $Z_1= \pi_G^{-1}(1) \cap Z$ is finite (because $Z$ is in $\lambda$) and normal in $Z$. Indeed, it is the kernel of the group homomorphism $Z \rightarrow G$. Then, let $Z$ act on $Z_1$ by conjugation. Since $Z$ is connected and $Z_1$ is finite, this action is trivial, i.e. $Z_1$ is central in $Z$. Let $H_0:= \pi_H(Z_1)$. So, $H_0$ is a finite subgroup of $H$. Let $H_1 \leq H$ be the centralizer in $H$ of $H_0$. Then, since $Z$ centralizes $Z_1$, we know that $\pi_H(Z)$ centralizes $\pi_H(Z_1)=H_0$, i.e. $\pi_H(Z) \leq H_1$. In other words, we have $Z \leq G \times H_1$.
			

			Now, let $f : G \times H_1 \rightarrow G \times (H_1/ H_0)$ be  the quotient morphism, and let $T=f(Z) \leq G \times (H_1 / H_0)$.
			
			Note that $\pi_G(T) = \pi_G(Z)$ is not in $\mu_G$, so by \cite[Lemma 3.17]{Stabilizers-NTP2}, it contains the connected component $G^{00}_{N_1}$ of $G$. Also, by construction, the fiber over $1_G$ in $T$ is trivial, so $T$ induces the graph of a group homomorphism $g: G^{00}_{N_1} \rightarrow H_1 / H_0$.

			Finally, as $Z$ is in $\lambda$, its fibers over elements of $H$ are in $\nu$. Then, since $\pi_H(Z_1)$ is finite, the fiber in $T$ over $1_{H^{\prime}}$ is a finite union of sets in $\nu$, so is in $\nu$ as well. In other words, the kernel of the morphism $g$ is in $\nu$, as required.
		\end{proof}

	\end{subsection}
	


	\begin{subsection}{Group configurations}

		Here, we describe another useful tool for identifying definable groups : group configurations. We say that two families $(a_i)_{i \in I}$ and $(b_i)_{i \in I}$ of tuples are \emph{equivalent} over a parameter set $A$ if, for all $i \in I$, we have $acl(Aa_i) = acl(Ab_i)$.

		In this section, we prove a general comparison theorem (Theorem \ref{theo_group_with_kernel_blind_to_S}) for groups definable in fields, inspired by Theorem \ref{theo_embedding_to_H(M)} below, but allowing for imaginaries. The main tools are the stabilizer theorem, which we use in Theorem \ref{theo_how_to_find_group_homomorphisms}, and the (variant of the) group configuration Theorem \ref{theo_relative_stable_group_config}.
		First, we need to introduce a natural condition on the algebraic closure operator.

		\begin{defi}\label{defi_alg_bounded} Let $\mathcal{L}$ be a language extending the language of rings. Let $T$ be an $\mathcal{L}$-theory, extending the theory of fields. We say that $T$ is \emph{algebraically bounded} if the following property holds : 
			

			For all $M \models T$, for all parameter sets $A \subseteq K(M)$, the algebraic closure $acl^M(A) \cap K$ coincides with the field-theoretic relative algebraic closure $K(A)^{alg} \cap K(M)$.
			
			
		\end{defi}

		\begin{theo}[\cite{Stabilizers-NTP2}, Theorem 2.19]\label{theo_embedding_to_H(M)} Let $T$ be an algebraically bounded theory of (possibly enriched) perfect fields. Let $G$ be a group definable in a sufficiently saturated model $M$ of $T$. Assume that $T$ admits an $M$-invariant ideal $\mu_G$ on $G$, stable under left and right translations, 
			which is $S1$ on $G$. Let $p \in S_G(M)$ be a type which is not in $\mu_G$, such that condition (F) holds : there are $a$, $b$ realizing $p$ such that $b \downfree_M a$ and $a \downfree_M b$.

			Then, there is an algebraic group $H$ over $M$, and a definable finite-to-one group homomorphism from a type-definable subgroup $D$ of $G$ to the definable group $H(M)$, such that $D$ is not in $\mu_G$.  
		\end{theo}

		\begin{defi}
			Let $A$ be a set of parameters.   A \emph{regular group configuration} over $A$ is a tuple $(\alpha_1, \alpha_2, \alpha_3, \beta_1, \beta_2, \beta_3)$ of elements satisfying the following properties :

			
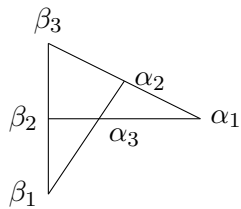
\begin{center}
				\begin{tikzpicture}
					\coordinate (b3) at (0,0) ;
					\coordinate (b2) at (0,-1) ;
					\coordinate (b1) at (0,-2) ;
					\coordinate (a2) at (1,-0.5) ;
					\coordinate (a1) at (2,-1) ;
					\coordinate (a3) at (1,-1) ;

					\draw (b3) -- (a1);
					\draw (b3) -- (b1);
					\draw (b2) -- (a1);
					\draw (b1) -- (a2);
					
					
					\draw (b3) node [above] {\(\beta_{3}\)};
					\draw (b2) node [left] {\(\beta_{2}\)};
					\draw (b1) node [left] {\(\beta_{1}\)};
					\draw (a2) node [right] {\(\alpha_{2}\)};
					\draw (a1) node [right] {\(\alpha_{1}\)};
					\draw (a3) node [below] {\(\alpha_{3}\)};
					
				\end{tikzpicture}
			\end{center}
			
			\begin{enumerate}    
				\item If $a,b,c$ are three non-colinear points in the diagram above, then the triple $(a,b,c)$ is an independent family over $A$, in the sense of (non)forking.
				\item If $a,b,c$  are three colinear points in the diagram above, then $a \in \mathrm{acl}(Abc)$.
				
			\end{enumerate}
			
			%
			%
			
		\end{defi}

		In this paper, we will only be interested in group configurations when the ambient theory is stable. Usually, this will mean considering a $6$-tuple of elements and a stable theory $T_0$, such that this $6$-tuple is a group configuration in the sense of $T_0$.

		\begin{prop}\label{prop_purely_im_sets}
			Let $T$ be an $\mathcal{L}$-theory, and $\mathcal{S}$ be a collection of $\mathcal{L}$-definable sets, closed under finite products. Then, 
			\begin{enumerate}
				
				\item The collection of definable sets (possibly with parameters in models of $T$) that are blind to $\mathcal{S}$ is closed under finite products.
				
				\item Assume that finite subsets of sets in $\mathcal{S}$ can be always be coded in $\mathcal{S}$. Then, a definable set $X$ is blind to $S$ if and only if, for all $S$ in $\mathcal{S}$, all definable functions from $X$ to $S$ have finite image.
				
			\end{enumerate}

		\end{prop}

		\begin{proof}
			Let us prove the first point. Let $X,Y$ be definable sets blind to $\mathcal{S}$. Let $S \in \mathcal{S}$. Let $F : X \times Y \rightarrow S$ be a definable finite correspondence. Then, for all $x \in X$, we get a definable finite correspondence $F_x : Y \rightarrow S$. Since $Y$ is blind to $\mathcal{S}$, this correspondence has a finite image $S_x$. Then, one can check that the set $\lbrace (x, s) \, | \, x \in X, s \in S_x \rbrace \subseteq X \times S$ is a definable finite correspondence from $X$ to $S$. Since $X$ is blind to $\mathcal{S}$, this correspondence has a finite image $S_0 \subseteq S$. Then, by construction, the image of $F$ is contained in $S_0$, which is finite, as required.
			
			For the second point, it suffices to note that, if $F: X \rightarrow S$ is a finite correspondence, and $k < \omega$ is such that fibers $F_x$ have size at most $k$, then, by coding subsets of $S$ of size at most $k$ in some $T \in \mathcal{S}$, one can find a definable function $f: X \rightarrow T$, whose image is infinite if and only if the image of $F$ is. This concludes the proof. 
		\end{proof}

		The first ingredient for the upcoming Theorem \ref{theo_relative_stable_group_config} is the construction, from a definable group, of a group configuration in a superstable quantifier-free reduct.

		\begin{lemma}\label{lemma_preserve_nonforking}
			
			Let $T_0, T_1$ be theories, in languages $\mathcal{L}_0 \subseteq \mathcal{L}_1$ respectively. Assume that $T_0$ is stable, has quantifier elimination and elimination of imaginaries (in $\mathcal{L}_0$). Also assume that $T_1|_{\mathcal{L}_0} \supseteq {T_0}_{\forall}$. Let $N_0 \models T_0$, and $N_1 \leq N_0$ be a model of $T_1$. In other words, $N_1$ is a model of $T_1$, and we are also given an $\mathcal{L}_0$-embedding $N_1 \rightarrow N_0$.

			Let $a,b,C$ be subsets of $N_1$. Assume that, in the sense of $T_1$, we have $a \downfree_C b$, and that $C=acl_1(C)$ is an extension base (for nonforking) in $T_1$. Then, in $T_0$, we have $a \downfree^0_C b$, i.e. $tp_0(a / Cb) = qftp_0(a/Cb)$ does not fork over $C$.

		\end{lemma}
		
		\begin{proof}
			The proof is copied from that of \cite[Lemme 2.1]{BPW-GeomRel}. First, up to taking an elementary extension of the pair $(N_0, N_1)$, we may assume that the models are sufficiently saturated. Now, let $\phi(x,y)$ be an $\mathcal{L}_0(C)$-formula, which we assume to be quantifier-free, satisfied by $(a,b)$. We wish to show that $\phi(x, b)$ does not fork, in the sense of $T_0$, over $C$. By stability (simplicity would be enough), it suffices to show that, for some Morley sequence $(b_i)_{i < \omega}$ in $tp_0(b / C)$, the collection of formulas $\lbrace \phi(x, b_i) \, | \, i < \omega \rbrace$ is consistent (in $T_0$).

			Since $C$ is an extension base, we may find, in $N_1$, large sequences $(b_i)_{i < \kappa}$ of realizations of $tp(b / C)$, such that $b_i \downfree_C b_{< i}$ for all $i < \kappa$. Then, by left transitivity for nonforking, we have $(b_j)_{ i \leq j < \kappa} \downfree_C (b_k)_{k < i}$ for all $i < \kappa$. Using Erdös-Rado, we can now extract a $C$-indiscernible subsequence $(b_i)_{i < \omega \cdot 2}$, which also satisfies these nonforking properties.

			On the one hand, as $(b_i)_{i < \omega \cdot 2}$ is $C$-indiscernible in $T_1$, with $b_0 \equiv_C b$ and $a \downfree_C b$, we know that $\lbrace \phi(x, b_i) \, | \, i < \omega \cdot 2 \rbrace$ is consistent in $T_1$. Then, this is also consistent in $T_0$, because $\phi$ is quantifier-free. Also, it is clear that $(b_i)_{i < \omega \cdot 2}$ is $C$-indiscernible in $T_0$. Then, the only thing left to prove is that $(b_i)_{i < \omega \cdot 2}$  is independent over $C$, in the sense of $T_0$.

			By $C$-indiscernibility and stability, the sequence $(b_i)_{ \omega \leq i < \omega \cdot 2}$ is a Morley sequence in the type $tp_0(b_{\omega} / C (b_j)_{j < \omega})$. In particular, it is $T_0$-independent over $C (b_j)_{j < \omega}$. Also, the canonical basis of the type $tp_0(b_{\omega} / C (b_j)_{j < \omega})$, which exists by elimination of imaginaries in $T_0$, is algebraic over $C (b_j)_{j < \omega}$, and definable over $C (b_i)_{ \omega \leq i < \omega \cdot 2}$. In particular, it belongs to the intersection
			
			\noindent $acl_1(C (b_j)_{j < \omega}) \cap acl_1(C (b_i)_{ \omega \leq i < \omega \cdot 2})$. However, the latter set is equal to $acl_1(C) = C$, because we have $(b_i)_{ \omega \leq i < \omega \cdot 2} \downfree_C (b_j)_{j < \omega}$ (inside $T_1$). This implies that $(b_i)_{ \omega \leq i < \omega \cdot 2}$ is in fact a Morley sequence in $tp_0(b_{\omega} / C)$, which concludes the proof. \end{proof}

		\begin{prop}\label{prop_how_to_build_group_configs_in_stable_reducts}
			
			Let $T_0, T_1$ be $\mathcal{L}_0 \subseteq \mathcal{L_1}$-theories, and $M_1 \preceq^+ N_1 \models T_1$. Let $S$ be an $\mathcal{L}_0$-quantifier-free-definable set.

			Assume that $T_0$ is superstable, and has elimination of quantifiers and imaginaries. Also assume that $T_1|_{\mathcal{L}_0} \supseteq {T_0}_{\forall}$. Let $N_0$ be a model of $T_0$, and $\iota : N_1 \rightarrow N_0$ be an $\mathcal{L}_0$-embedding, such that $dcl_0(N_1) = N_1$ and $acl_0(A) \cap N_1 = acl_1(A)$ for all $A \subseteq N_1$.

			Let $G$ be an $M_1$-definable group in $T_1$, with an $M_1$-definable map $\pi: G \rightarrow S$ whose fibers are blind to $S$. Assume that there exists a type $p \in S_G(N_1)$ which is strongly right-f-generic over $M_1$. Let $a \models p|_{M_1}$, $b \models p|_{M_1a}$, $c \models p|_{M_1ab}$. Then, the following tuple is a regular group configuration over $M_1$, in the sense of $T_0$ : $(\pi(a\cdot c),\pi(c), \pi(b \cdot a \cdot c),\pi(b \cdot a), \pi(b), \pi(a))$.

		\end{prop}
		
		\begin{proof}
			Let $\downfree^0$ denote independence in the superstable reduct given by quantifier-free $\mathcal{L}_0$-formulas, over $M_1$. To simplify notations, let $\alpha$, $\beta$, $\gamma$, $\alpha \cdot \gamma$, $\beta \cdot \alpha$, $\beta\cdot \alpha \cdot \gamma$ denote $\pi(a)$, $\pi(b)$, $\pi(c)$, $\pi(a\cdot c)$, $\pi(b \cdot a)$, $\pi(b \cdot a \cdot c)$ respectively.
			
			We wish to show that the following is a regular group configuration over $M_1$, in the sense of $N_0$:

			
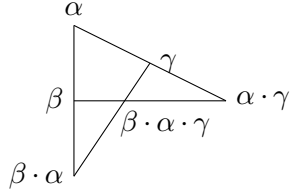
\begin{center}
				\begin{tikzpicture}
					
					\coordinate (b3) at (0,0) ;
					\coordinate (b2) at (0,-1) ;
					\coordinate (b1) at (0,-2) ;
					\coordinate (a2) at (1,-0.5) ;
					\coordinate (a1) at (2,-1) ;
					\coordinate (a3) at (1.2,-1) ;

					\draw (b3) -- (a1);
					\draw (b3) -- (b1);
					\draw (b2) -- (a1);
					\draw (b1) -- (a2);
					
					
					\draw (b3) node [above] {\(\alpha\)};
					\draw (b2) node [left] {\(\beta\)};
					\draw (b1) node [left] {\(\beta \cdot \alpha\)};
					\draw (a2) node [right] {\(\gamma\)};
					\draw (a1) node [right] {\(\alpha \cdot \gamma\)};
					\draw (a3) node [below] {\(\beta\cdot \alpha \cdot \gamma\)};

				\end{tikzpicture}
			\end{center}

			\begin{claim}
				Any three colinear points in the diagram above are interalgebraic over $M_1$ in $T_1$. 
			\end{claim}
			
			\begin{proof}
				We use the hypothesis that the fibers of $\pi$ are blind to $S$. For instance, we know that $b\cdot a \in acl_1(M_1, b, a)$, and we wish to deduce that $\pi(b \cdot a) \in acl_1(M_1, \pi(b), \pi(a))$. Consider the following $M_1 \pi(a) \pi(b)$-definable set: $F:= \lbrace (g_1, g_2, s) \in G \times G \times S \, | \, \pi(g_1) = \pi(a) \, \wedge \pi(g_2) = \pi(b) \, \wedge \pi(g_2 \cdot g_1) = s \rbrace$. Then, $F$ induces a partial definable map $S_{\pi(a)} \times S_{\pi(b)} \rightarrow S$, whose image contains $\pi(b \cdot a)$. Since the fibers $S_{\pi(a)}$ and $S_{\pi(b)}$ are blind to $S$, the image of $F$ is finite. Thus, we have found a finite $M_1 \pi(a) \pi(b)$-definable set containing $\pi(b \cdot a)$, as required.
			\end{proof}
			
			Then, using the hypothesis on $acl_0$, any three colinear points are $T_0$-interalgebraic over $M_1$. So, it remains to prove that the required $\downfree^0$-independence relations hold. To do so, we will use $U$-rank computations, performed in $N_0$.

			Note that, thanks to Lemma \ref{lemma_preserve_nonforking}, we can deduce several $\downfree^0$-independence relations from the construction.
			Let $\lambda$ denote the ordinal which is the $U$-rank of $\alpha$ over $M_1$, computed in the superstable structure $N_0$. Note that $\lambda = U(\alpha / M_1) = U(\beta / M_1) = U(\gamma / M_1)$. This is because $tp(a / M_1) = tp(b/M_1) = tp(c/M_1) = p$ and so $tp_0(\pi(a) / M_1) = tp_0(\pi(b)/M_1) = tp_0(\pi(c)/M_1)$.

			\begin{claim}
				
				The $U$-ranks over $M_1$ of the elements $\alpha \cdot \gamma, \beta \cdot \alpha, \beta\cdot \alpha \cdot \gamma$ are all equal to $\lambda$ as well.
				
			\end{claim}
			
			\begin{proof}
				By genericity and Lemma \ref{lemma_preserve_nonforking}, we know that $\beta \cdot \alpha \downfree^0 \alpha$. So, by Lascar's equality (see \cite[Theorem 19.5]{Poi-MT}), we know that $U(\beta \cdot \alpha, \alpha / M_1) = U(\beta\cdot \alpha / M_1) \oplus U(\alpha / M_1) = U(\beta\cdot \alpha / M_1) \oplus \lambda$. Then, by interalgebraicity, we deduce equality of the ranks : $U(\beta \cdot \alpha, \alpha / M_1) = U(\beta, \alpha / M_1) = \lambda \oplus \lambda$, using Lascar's equality again, since $\beta \downfree^0 \alpha$. As all ordinals are cancellable for the symmetric sum $\oplus$, we deduce that $U(\beta \cdot \alpha / M_1) = \lambda$.
				
				Similarly, using the facts $\alpha \cdot \gamma \downfree^0 \alpha$ and $\gamma \downfree^0 \alpha$, and the algebraicities, we compute $\lambda \oplus \lambda = U(\alpha, \gamma/M_1)= U(\alpha \cdot \gamma, \gamma / M_1) = U(\alpha \cdot \gamma / M_1) \oplus U(\gamma / M_1) = U(\alpha \cdot \gamma / M_1) \oplus \lambda$. So, as above, $U(\alpha \cdot \gamma / M_1) = \lambda$. Finally, as $\gamma \downfree^0 \beta \cdot \alpha$ and $\beta \cdot \alpha \cdot \gamma \downfree^0 \beta \cdot \alpha $, the same arguments show that $U(\beta \cdot \alpha \cdot \gamma / M_1) = U(\beta / M_1) = \lambda$.
			\end{proof}

			Now, let $x,y,z$ be non-collinear elements in the configuration (for instance, $x= \alpha$, $y=\alpha \cdot \gamma$, and $z=\beta \cdot \alpha$). Then, we know that $acl_0(M_1 x y z) = acl_0(M_1 \alpha  \beta  \gamma)$, so $U(x, y, z / M_1) = U(\alpha, \beta, \gamma / M_1) = \lambda \oplus \lambda \oplus \lambda$. Then, using Lascar's inequality, we compute:

			\begin{align*}
				U(x, y, z / M_1) &\leq U(x / M_1 y  z) \oplus U(y / M_1 z) \oplus U(z / M_1) \\
				&\leq U(x / M_1) \oplus U(y / M_1) \oplus U(z / M_1) \\
				&= \lambda \oplus \lambda \oplus \lambda \\
				&= U(x, y, z / M_1).
			\end{align*}

			So, we have the equalities $U(x / M_1 y  z) = U(x / M_1)$ and $U(y / M_1 z) = U(y / M_1)$, which precisely mean that $x, y, z$ is an independent family over $M_1$, as required.
		\end{proof}

		\begin{lemma}\label{lemma_fin_sat_implies_stationary}

			Let $T_0$ be a stable theory. Let $A \subseteq M \models T_0$, and $p \in S(A)$. Assume that $p$ is finitely satisfiable in $A$. Then it is stationary. Moreover, for all $B \supseteq A$, the type $p|_B$ is finitely satisfiable in $A$.
			
		\end{lemma}
		
		\begin{proof}
			We may assume that $M$ is $|A|^+$-saturated and $|A|^+$-strongly homogeneous. Using compactness, it is straightforward to check that $p$ admits an extension $q \in S(M)$ which is also finitely satisfiable in $A$. Then, $q$ is $A$-invariant, so $A$-definable, and $q|_A = p$ is stationary.

			Moreover, if $B \supseteq A$, then some nonforking extension of $p$ to $B$ is finitely satisfiable in $A$; we conclude by uniqueness. \end{proof}

		The following theorem was inspired by \cite[Proposition 3.1]{HruPil-GpPFF}.

		\begin{theo}\label{theo_relative_stable_group_config}
			
			Let $T_0, T_1$ be theories, in languages $\mathcal{L}_0 \subseteq \mathcal{L}_1$ respectively. Assume that $T_0$ is superstable, has quantifier elimination and elimination of imaginaries in $\mathcal{L}_0$, in a collection of sorts $\mathcal{S}$.  Assume that $T_1$ is $\mathrm{NIP}$. Also assume that $T_1|_{\mathcal{L}_0} \supseteq {T_0}_{\forall}$. Let $N_0 \models T_0$, ${N_1} \models T_1$ be very saturated, and $\iota : N_1 \rightarrow N_0$ be an $\mathcal{L}_0$-embedding. Assume that $dcl_0(N_1) = N_1$ and $acl_0(A) \cap N_1 = acl_1(A)$ for all $A \subseteq N_1$.
			Let $V$ be a product of sorts of $\mathcal{S}$. Let $G$ be a group definable in ${N_1}$, with a definable map $\pi: G \rightarrow V$, whose fibers are blind to $\mathcal{S}$. Let ${M_1} \preceq^+ {N_1}$, with $M_1$ sufficiently saturated, and let $p \in S_G(M_1)$ be strongly bi-f-generic over some small model $C \preceq^+ M_1$. Let ${a} \models p|_{M_1}$, ${b} \models p|_{M_1 {a}}$ and ${c} \models p|_{M_1 {a}{b}}$.

			Then, there exists an $\mathcal{L}_0(M_1)$-quantifier-free-type-definable group $\Gamma$, and elements $g_1 \in \Gamma(N_1)$, $g_2 \in \Gamma(N_1)$, $g_3 \in \Gamma(N_1)$ that are independent and generic over $M_1$, such that the following configurations are $acl_1$-equivalent over $M_1$:

			\begin{center}
				\begin{tikzpicture}
					
					\coordinate (b3) at (0,0) ;
					\coordinate (b2) at (0,-1) ;
					\coordinate (b1) at (0,-2) ;
					\coordinate (a2) at (1,-0.5) ;
					\coordinate (a1) at (2,-1) ;
					\coordinate (a3) at (1.2,-1) ;

					\draw (b3) -- (a1);
					\draw (b3) -- (b1);
					\draw (b2) -- (a1);
					\draw (b1) -- (a2);

					\coordinate (b3+1) at (5,0) ;
					\coordinate (b2+1) at (5,-1) ;
					\coordinate (b1+1) at (5,-2) ;
					\coordinate (a2+1) at (6,-0.5) ;
					\coordinate (a1+1) at (7,-1) ;
					\coordinate (a3+1) at (6.3,-1) ;

					\draw (b3+1) -- (a1+1);
					\draw (b3+1) -- (b1+1);
					\draw (b2+1) -- (a1+1);
					\draw (b1+1) -- (a2+1);

					\draw (b3) node [above] {\(\pi(a)\)};
					\draw (b2) node [left] {\(\pi(b)\)};
					\draw (b1) node [left] {\(\pi(b \cdot a)\)};
					\draw (a2) node [right] {\(\pi(c)\)};
					\draw (a1) node [right] {\(\pi(a \cdot c)\)};
					\draw (a3) node [below] {\(\pi(b\cdot a \cdot c)\)};

					\draw (b3+1) node [above] {\(g_1\)};
					\draw (b2+1) node [left] {\(g_2\)};
					\draw (b1+1) node [left] {\(g_2 \cdot g_1\)};
					\draw (a2+1) node [right] {\(g_3\)};
					\draw (a1+1) node [right] {\(g_1 \cdot g_3\)};
					\draw (a3+1) node [below] {\(g_2\cdot g_1 \cdot g_3\)};

				\end{tikzpicture}
			\end{center}

		\end{theo}

		\begin{rem}\label{rem_G(N)_is_a_group}
			
			\begin{enumerate}
				
				\item Since $T_1$ is NIP, the existence of a strong bi-f-generic $p$ is implied by, and in fact turns out to be equivalent to, definable amenability of $G$. See Proposition \ref{prop_exists_generics}.
				
				\item In the conclusion of the theorem, since $N_1$ is a definably closed $\mathcal{L}_0$-substructure of $N_0$, and $\Gamma_0$ is $\mathcal{L}_0(M_1)$-quantifier-free-definable, the set $\Gamma_0(N_1) \subseteq \Gamma_0(N_0)$ is in fact a subgroup, and thus $\Gamma_0(N_1)$ is an $\mathcal{L}_0(M_1)$-quantifier-free-definable group, in $N_1$ itself.

				\item The hypotheses of the theorem imply that $M_1$ is also $dcl_0$-closed in $N_0$.  
				
			\end{enumerate}
			
		\end{rem}

		\begin{proof}[Proof of Theorem \ref{theo_relative_stable_group_config}]
			First, since working with projections can be tedious, we shall change our notations slightly. Let us write ${C_{\mathcal{S}}}, {M_{1, {\mathcal{S}}}}, {N_{1, {\mathcal{S}}}}$ instead of $C|_{\mathcal{L}_0}, M_1|_{\mathcal{L}_0}, N_1|_{\mathcal{L}_0}$. Then, if ${d}$, ${B}$, etc. denote tuples or parameters which are in ${N_1}$, we let $d_{\mathcal{S}}$, $B_{\mathcal{S}}$, etc. denote their projections to the sorts in $\mathcal{S}$. Abusing notations slightly, if ${\alpha}, {\beta}$ are elements in $G({N_1})$, we write $\alpha_{\mathcal{S}} \cdot  \beta_{\mathcal{S}}$ to denote the image $\pi({\alpha} \cdot {\beta}) \in V$. We shall keep these conventions until the end of the proof of the theorem.

			The main point of this proof is to follow the recipe given by the proof of the stable group configuration theorem, while making sure that the elements and parameters stay inside $N_{1, {\mathcal{S}}}$. \emph{However, trying to apply the stable group configuration theorem directly does not seem to work.} Instead, most of the work will be done inside the $\mathrm{NIP}$ theory $T_1$, exploiting the strongly bi-f-generic type $p$ 
			several times.  It is crucial to keep in mind that there are two theories at play here. As in the statement of the theorem, we let $dcl_0$, $acl_0$, $tp_0$, etc. denote the $\mathcal{L}_0$-quantifier-free notions, i.e. those defined by the superstable theory $T_0$.

			Let ${b_{1}} = {b} \cdot {a}$, ${b_2} = {b}$, ${b_3} = {a}$, ${a_1} = {a} \cdot {c}$, ${a_2} = {c}$ and ${a_3} = {b} \cdot {a}\cdot {c}$. By Proposition \ref{prop_how_to_build_group_configs_in_stable_reducts}, the following is, in $T_0$, a regular group configuration over both $M_{1, \mathcal{S}}$ and $C_{\mathcal{S}}$:

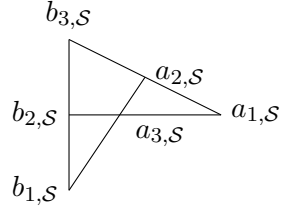
			\begin{center}
				\begin{tikzpicture}
					
					\coordinate (b3) at (0,0) ;
					\coordinate (b2) at (0,-1) ;
					\coordinate (b1) at (0,-2) ;
					\coordinate (a2) at (1,-0.5) ;
					\coordinate (a1) at (2,-1) ;
					\coordinate (a3) at (1.2,-1) ;

					\draw (b3) -- (a1);
					\draw (b3) -- (b1);
					\draw (b2) -- (a1);
					\draw (b1) -- (a2);
					
					
					\draw (b3) node [above] {\(b_{3,  \mathcal{S}}\)};
					\draw (b2) node [left] {\(b_{2,  \mathcal{S}}\)};
					\draw (b1) node [left] {\(b_{1,  \mathcal{S}}\)};
					\draw (a2) node [right] {\(a_{2,  \mathcal{S}}\)};
					\draw (a1) node [right] {\(a_{1,  \mathcal{S}}\)};
					\draw (a3) node [below] {\(a_{3,  \mathcal{S}}\)};

				\end{tikzpicture}
			\end{center}
			
			%
			%
			%
			%
			%
			%
			%
			%
			%

			%
				%
				%
				%
				%
				%
				%
				%
				%
				%
				%

			\begin{claim}\label{claim_wma_interdcl}
				There exists a small model ${C_2} \models T_1$, with ${C} \preceq {C_2} \preceq^+ {M_1}$, and a regular group configuration over $C_{2, \mathcal{S}}$ (in the sense of $T_0$), which is made of elements $(c_1, c_2, c_3, d_1, d_2, d_3)$  of $N_{1, \mathcal{S}}$, and equivalent over $C_{2, \mathcal{S}}$ to $(a_{1, \mathcal{S}}, a_{2, \mathcal{S}}, a_{3, \mathcal{S}}, b_{1, \mathcal{S}}, b_{2, \mathcal{S}}, b_{3, \mathcal{S}})$, with the following properties :
				
				\begin{enumerate}
					\item We have $c_1 \in dcl_0(C_{2, \mathcal{S}} d_2 c_3)$ and $c_2 \in dcl_0(C_{2, \mathcal{S}} d_1 c_3)$.
					
					\item The element $c_3$ is in $dcl_0(C_{2, \mathcal{S}} d_1 c_2) \cap dcl_0(C_{2, \mathcal{S}} d_2 c_1)$.
					
				\end{enumerate}

			\end{claim}
			
			\begin{proof}[Proof of the claim]
				We begin with the first definability property. Let ${b^{\prime}_1}$ and  ${b^{\prime\prime}_2}$ be elements of $G({M_1})$ realizing $p|_{{C}}$. Then, let ${b^{\prime}_3} = {b^{}_2}^{-1} \cdot {b^{\prime}_1}$, ${a^{\prime}_2}= {b^{\prime}_1}^{-1} \cdot {a^{}_3}$, ${b^{\prime\prime}_3} = {b^{\prime\prime}_2}^{-1} \cdot {b^{}_1}$, and ${a^{\prime\prime}_1}={b^{\prime\prime}_2}^{-1} \cdot {a^{}_3}$. In other words, we apply Proposition \ref{prop_how_to_build_group_configs_in_stable_reducts}, up to a permutation of the elements, to get the following regular group configurations over $C_{\mathcal{S}}$, and also over $M_{1, \mathcal{S}}$ :

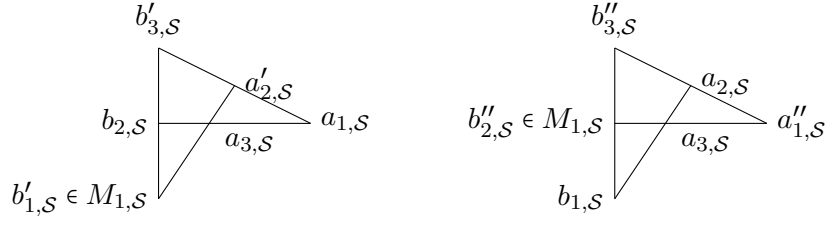
				\begin{center}
					\begin{tikzpicture}

						\coordinate (b3) at (0,0) ;
						\coordinate (b2) at (0,-1) ;
						\coordinate (b1) at (0,-2) ;
						\coordinate (a2) at (1,-0.5) ;
						\coordinate (a1) at (2,-1) ;
						\coordinate (a3) at (1.2,-1) ;

						\draw (b3) -- (a1);
						\draw (b3) -- (b1);
						\draw (b2) -- (a1);
						\draw (b1) -- (a2);

						\coordinate (b3+1) at (6,0) ;
						\coordinate (b2+1) at (6,-1) ;
						\coordinate (b1+1) at (6,-2) ;
						\coordinate (a2+1) at (7,-0.5) ;
						\coordinate (a1+1) at (8,-1) ;
						\coordinate (a3+1) at (7.2,-1) ;

						\draw (b3+1) -- (a1+1);
						\draw (b3+1) -- (b1+1);
						\draw (b2+1) -- (a1+1);
						\draw (b1+1) -- (a2+1);

						\draw (b3) node [above] {\(   b^{\prime}_{3, \mathcal{S}}   \)};
						\draw (b2) node [left] {\(    b_{2, \mathcal{S}}  \)};
						\draw (b1) node [left] {\(    b^{\prime}_{1, \mathcal{S}} \in M_{1, \mathcal{S}}  \)};
						\draw (a2) node [right] {\(   a^{\prime}_{2 , \mathcal{S}}   \)};
						\draw (a1) node [right] {\(   a_{1, \mathcal{S}}  \)};
						\draw (a3) node [below] {\(   a_{3, \mathcal{S}} \)};

						\draw (b3+1) node [above] {\(   b^{\prime\prime}_{3, \mathcal{S}}   \)};
						\draw (b2+1) node [left] {\(    b^{\prime\prime}_{2, \mathcal{S}} \in M_{1, \mathcal{S}}   \)};
						\draw (b1+1) node [left] {\(   b_{1, \mathcal{S}}   \)};
						\draw (a2+1) node [right] {\(   a^{}_{2, \mathcal{S}}   \)};
						\draw (a1+1) node [right] {\(   a^{\prime\prime}_{1, \mathcal{S}} \)};
						\draw (a3+1) node [below] {\(   a_{3, \mathcal{S}}    \)};
						
					\end{tikzpicture}
				\end{center}
				%
				%
					%
					%
					%
					%
					%
					%
					%
					%
					%
					%
					%
					%
					%
					%
					%
					%
					%
					%
					%
					%
					%
					%
					%
					%
				%
				%

				Then, let $\widetilde{a_1}$, resp. $\widetilde{a_2}$, denote the code of the set of $acl_0$-conjugates of $a_{1, \mathcal{S}}$ over $C  b^{\prime}_{3, \mathcal{S}}  a^{\prime}_{2, \mathcal{S}}  b^{}_{2, \mathcal{S}} a_{3, \mathcal{S}}$, resp. of $a_{2, \mathcal{S}}$ over $C  b^{\prime\prime}_{2, \mathcal{S}}  a^{\prime\prime}_{1, \mathcal{S}} b_{1, \mathcal{S}} a_{2, \mathcal{S}}$.     Note that, by elimination of imaginaries and $dcl_0$-closure of $N_{1, \mathcal{S}}$, we may assume that these elements are in $N_{1, \mathcal{S}}$. Also, since $a_{1, \mathcal{S}}$ belongs to the finite set coded by $\widetilde{a_1}$, we have $a_{1, \mathcal{S}} \in acl_0(C \widetilde{a_1})$. Similarly, we have $a_{2, \mathcal{S}} \in acl_0(C \widetilde{a_2})$. On the other hand, the algebraicities in $T_0$ imply that $\widetilde{a_1} \in acl_0(Cb_{2, \mathcal{S}} a_{3, \mathcal{S}}) \cap acl_0(C b^{\prime}_{3, \mathcal{S}} a^{\prime}_{2, \mathcal{S}})$, and $b_{2, \mathcal{S}} a_{3, \mathcal{S}} \downfree^0_{Ca_{1, \mathcal{S}}}  b^{\prime}_{3, \mathcal{S}} a^{\prime}_{2, \mathcal{S}}$. Hence, the element  $\widetilde{a_1}$ is in $acl_0(C a_{1, \mathcal{S}})$. Similarly, we have $acl_0(C \widetilde{a_2}) = acl_0(C a_{2, \mathcal{S}})$.

				What we gained in the process is the following: 
				\begin{center}
					$\widetilde{a_1} \in dcl_0(C b_{2, \mathcal{S}} a_{3, \mathcal{S}} b^{\prime}_{3, \mathcal{S}} a^{\prime}_{2, \mathcal{S}})$ and $\widetilde{a_2} \in dcl_0(C b_{1, \mathcal{S}} a_{3, \mathcal{S}} b^{\prime\prime}_{3, \mathcal{S}} a^{\prime\prime}_{1, \mathcal{S}})$
				\end{center}


				So, let ${C_1}$ be a small model of $T_1$ containing ${C} {b^{\prime}_1} {b^{\prime\prime}_2}$ such that ${C_1} \preceq^+ {M_1}$. Let $(\alpha_1, \alpha_2, \alpha_3, \beta_1, \beta_2, \beta_3) \in N_{1, \mathcal{S}}$ denote the following configuration :

				
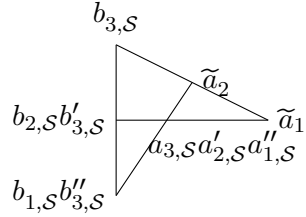
\begin{center}
					\begin{tikzpicture}
						
						\coordinate (b3) at (0,0) ;
						\coordinate (b2) at (0,-1) ;
						\coordinate (b1) at (0,-2) ;
						\coordinate (a2) at (1,-0.5) ;
						\coordinate (a1) at (2,-1) ;
						\coordinate (a3) at (1.4,-1) ;

						\draw (b3) -- (a1);
						\draw (b3) -- (b1);
						\draw (b2) -- (a1);
						\draw (b1) -- (a2);
						
						\draw (b3) node [above] {\(   b^{}_{3, \mathcal{S}}   \)};
						\draw (b2) node [left] {\(    b_{2, \mathcal{S}}b^{\prime}_{3, \mathcal{S}}  \)};
						\draw (b1) node [left] {\(    b^{}_{1, \mathcal{S}} b^{\prime\prime}_{3, \mathcal{S}}  \)};
						\draw (a2) node [right] {\(   \widetilde{a}^{}_{2}  \)};
						\draw (a1) node [right] {\(   \widetilde{a}_{1}  \)};
						\draw (a3) node [below] {\(   a_{3, \mathcal{S}} a^{\prime}_{2, \mathcal{S}} a^{\prime\prime}_{1, \mathcal{S}} \)};

					\end{tikzpicture}
				\end{center}
				\normalsize
				
				%
					%
					%
					%
					%
					%
					%
					%
					%
					%
					%
				
				This tuple is $acl_0$-equivalent to the previous one over $C_{1, \mathcal{S}}$, so it is indeed a regular group configuration over $C_{1, \mathcal{S}}$, and also over $M_{1, \mathcal{S}}$. Also, unfolding the definitions, we know that ${\beta_i} \in \mathcal{S}(acl_1({C_1} {b_i}))$ for $i=1,2,3$ and $ {\alpha_3} \in \mathcal{S}(acl_1({C_1} {a_3}))$. Regarding the $\beta_i$ and $b_i$, we even have definability over $C_1$, but we do not actually need it to finish the proof.

				Now, let us deal with the second definability requirement. Let $\widetilde{\alpha_3}$ denote a tuple in $N_{1, \mathcal{S}}$ which encodes the finite set of $acl_0$-conjugates of $\alpha_3$ over $C_{1, \mathcal{S}} \beta_1 \beta_2 \alpha_1 \alpha_2$. Then, for reasons similar to the previous step, we have $acl_0(C_{1, \mathcal{S}} \alpha_3) = acl_0(C_{1, \mathcal{S}} \widetilde{\alpha_3})$. Moreover, since $\alpha_1 \in dcl_0(C_{1, \mathcal{S}} \alpha_3 \beta_2)$ and  $\alpha_2 \in dcl_0(C_{1, \mathcal{S}} \alpha_3 \beta_1)$, and all elements in the finite set coded by $\widetilde{\alpha_3}$ have the same type as $\alpha_3$ over $C_{1, \mathcal{S}} \beta_1 \beta_2 \alpha_1 \alpha_2$, we still have $\alpha_1 \in dcl_0(C_{1, \mathcal{S}} \widetilde{\alpha_3} \beta_2)$ and  $\alpha_2 \in dcl_0(C_{1, \mathcal{S}}  \widetilde{\alpha_3} \beta_1)$.


				Recall that, by construction, the element $a_3$ is equal to $b_2\cdot b_3 \cdot a_2$, where $tp(a_2/M_1 b_1 b_2 b_3)$ is strongly bi-f-generic over $C_1$. Thus ${a_3} \downfree_{{C_1}} {M_1} {b_1} {b_2} {b_3}$.   So, by the algebraicities above, we have $ {\alpha_3}  \downfree_{{C_1}} {M_1} {\beta_1} {\beta_2} {\beta_3}$. Also, since $b_1 = b \cdot a$, $b_2 = b$ and $b_3 = a$, where $tp(b/ M_1 a)$ is strongly bi-f-generic over $C_1$, we know that ${b_1} \downfree_{{C_1}} M_1 {b_3}$ and ${b_2} \downfree_{{C_1}} M_1 {b_3}$. Since the $\beta_i$ are algebraic over the $b_i$, this implies ${\beta_1} \downfree_{{C_1}} M_1 {\beta_3}$ and ${\beta_2} \downfree_{{C_1}} M_1 {\beta_3}$.  Then, using left transitivity of nonforking (which always holds), we get ${\alpha_3} {\beta_1} \downfree_{{C_1}} {M_1} {\beta_3}$ and ${\alpha_3} {\beta_2} \downfree_{{C_1}} {M_1}{\beta_3}$.

				Now, by saturation of ${M_1}$, let ${\beta^{\prime\prime\prime}_3} \in {M_1}$ realize $tp({\beta_3} / {C_1})$. By $\mathrm{NIP}$, we know from the above independence relations that the types $tp({\alpha_3} {\beta_1} / {M_1} {\beta_3})$ and $tp({\alpha_3} {\beta_2} / {M_1} {\beta_3})$ admit a ${C_1}$-invariant extension (see for instance \cite[Corollary 5.29]{Sim-Book}). Thus, we have ${\alpha_3} {\beta_2}{\beta_3} \equiv_{{C_1}} {\alpha_3}{\beta_2}{\beta^{\prime\prime\prime}_3}$ and ${\alpha_3} {\beta_1}{\beta_3} \equiv_{{C_1}} {\alpha_3}{\beta_1}{\beta^{\prime\prime\prime}_3}$. Using the definability properties we just proved, we deduce ${\alpha_3} {\beta_2}{\beta_3} \alpha_1 \equiv_{{C_1}} {\alpha_3}{\beta_2}{\beta^{\prime\prime\prime}_3} \alpha_1$ and ${\alpha_3} {\beta_1}{\beta_3}\alpha_2 \equiv_{{C_1}} {\alpha_3}{\beta_1}{\beta^{\prime\prime\prime}_3} \alpha_2$. Then, by homogeneity of $N_1$, let ${\beta^{\prime\prime\prime}_1}$, ${\beta^{\prime\prime\prime}_2}$, $\alpha^{\prime\prime\prime}_1$ and $\alpha^{\prime\prime\prime}_2$ be elements in ${N_1}$ such that $ {\beta_1}{\beta_2}{\beta_3}\alpha_1 \alpha_2{\alpha_3} \equiv_{{C_1}} {\beta^{\prime\prime\prime}_1}{\beta_2}{\beta^{\prime\prime\prime}_3}\alpha_1 \alpha^{\prime\prime\prime}_2{\alpha_3} \equiv_{{C_1}} {\beta_1}{\beta^{\prime\prime\prime}_2}{\beta^{\prime\prime\prime}_3}\alpha^{\prime\prime\prime}_1 \alpha_2{\alpha_3}$.

				In particular, we have $ {\beta_1}{\beta_2}{\beta_3}\alpha_1 \alpha_2\widetilde{\alpha_3} \equiv^0_{C_{1, \mathcal{S}}} {\beta^{\prime\prime\prime}_1}{\beta_2}{\beta^{\prime\prime\prime}_3}\alpha_1 \alpha^{\prime\prime\prime}_2\widetilde{\alpha_3} \equiv^0_{C_{1, \mathcal{S}}} {\beta_1}{\beta^{\prime\prime\prime}_2}{\beta^{\prime\prime\prime}_3}\alpha^{\prime\prime\prime}_1 \alpha_2\widetilde{\alpha_3}$. In other words, we managed to build three copies of the configuration inside $N_{1, \mathcal{S}}$:

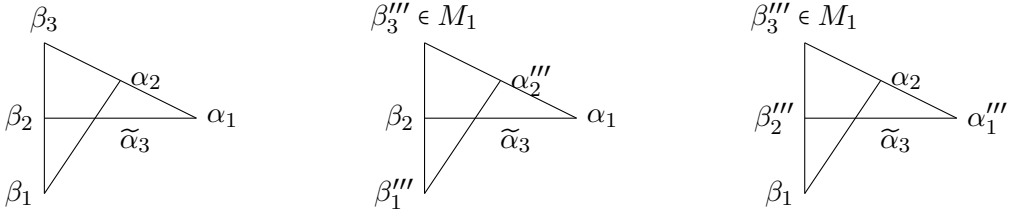
				\begin{center}
					\begin{tikzpicture}

						\coordinate (b3) at (0,0) ;
						\coordinate (b2) at (0,-1) ;
						\coordinate (b1) at (0,-2) ;
						\coordinate (a2) at (1,-0.5) ;
						\coordinate (a1) at (2,-1) ;
						\coordinate (a3) at (1.2,-1) ;

						\draw (b3) -- (a1);
						\draw (b3) -- (b1);
						\draw (b2) -- (a1);
						\draw (b1) -- (a2);

						\coordinate (b3+1) at (5,0) ;
						\coordinate (b2+1) at (5,-1) ;
						\coordinate (b1+1) at (5,-2) ;
						\coordinate (a2+1) at (6,-0.5) ;
						\coordinate (a1+1) at (7,-1) ;
						\coordinate (a3+1) at (6.2,-1) ;

						\draw (b3+1) -- (a1+1);
						\draw (b3+1) -- (b1+1);
						\draw (b2+1) -- (a1+1);
						\draw (b1+1) -- (a2+1);

						\coordinate (b3+2) at (10,0) ;
						\coordinate (b2+2) at (10,-1) ;
						\coordinate (b1+2) at (10,-2) ;
						\coordinate (a2+2) at (11,-0.5) ;
						\coordinate (a1+2) at (12,-1) ;
						\coordinate (a3+2) at (11.2,-1) ;

						\draw (b3+2) -- (a1+2);
						\draw (b3+2) -- (b1+2);
						\draw (b2+2) -- (a1+2);
						\draw (b1+2) -- (a2+2);

						\draw (b3) node [above] {\(\beta_3\)};
						\draw (b2) node [left] {\(\beta_2\)};
						\draw (b1) node [left] {\(\beta_1\)};
						\draw (a2) node [right] {\(\alpha_2\)};
						\draw (a1) node [right] {\(\alpha_1\)};
						\draw (a3) node [below] {\(\widetilde{\alpha}_3\)};

						\draw (b3+1) node [above] {\(\beta^{\prime\prime\prime}_{3} \in M_1\)};
						\draw (b2+1) node [left] {\(\beta_{2}\)};
						\draw (b1+1) node [left] {\(\beta^{\prime\prime\prime}_{1}\)};
						\draw (a2+1) node [right] {\(\alpha^{\prime\prime\prime}_{2}\)};
						\draw (a1+1) node [right] {\(\alpha_{1}\)};
						\draw (a3+1) node [below] {\(\widetilde{\alpha}_{3}\)};

						\draw (b3+2) node [above] {\(\beta^{\prime\prime\prime}_{3} \in M_1\)};
						\draw (b2+2) node [left] {\(\beta^{\prime\prime\prime}_{2}\)};
						\draw (b1+2) node [left] {\(\beta_{1}\)};
						\draw (a2+2) node [right] {\(\alpha_{2}\)};
						\draw (a1+2) node [right] {\(\alpha^{\prime\prime\prime}_{1}\)};
						\draw (a3+2) node [below] {\(\widetilde{\alpha}_{3}\)};
						
					\end{tikzpicture}
				\end{center}

				To conclude the proof of the claim, let ${C_2} \preceq^+ {M_1}$ be some small model of $T_1$ containing ${C_1} {\beta^{\prime\prime\prime}_3}$. We claim that the following configuration, over $C_{2, \mathcal{S}}$, has the required properties.

				
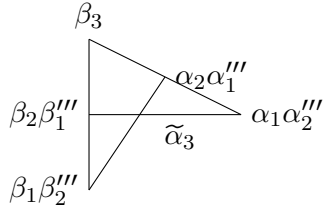
\begin{center}
					\begin{tikzpicture}
						
						\coordinate (b3) at (0,0) ;
						\coordinate (b2) at (0,-1) ;
						\coordinate (b1) at (0,-2) ;
						\coordinate (a2) at (1,-0.5) ;
						\coordinate (a1) at (2,-1) ;
						\coordinate (a3) at (1.2,-1) ;

						\draw (b3) -- (a1);
						\draw (b3) -- (b1);
						\draw (b2) -- (a1);
						\draw (b1) -- (a2);

						\draw (b3) node [above] {\(\beta_3\)};
						\draw (b2) node [left] {\(\beta_2\beta^{\prime\prime\prime}_{1}\)};
						\draw (b1) node [left] {\(\beta_1\beta^{\prime\prime\prime}_{2}\)};
						\draw (a2) node [right] {\(\alpha_2\alpha^{\prime\prime\prime}_{1}\)};
						\draw (a1) node [right] {\(\alpha_1\alpha^{\prime\prime\prime}_{2}\)};
						\draw (a3) node [below] {\(\widetilde{\alpha}_3\)};
						
					\end{tikzpicture}
				\end{center}
				%
				%
				%
					%
					%
					%
					%
					%
					%
					%
					%
					%
					%

				First, it is $acl_0$-equivalent over $C_2$, to the previous ones, so it is indeed a regular group configuration over $C_{2, \mathcal{S}}$, and also over $M_{1, \mathcal{S}}$. Then, by construction, we have $\widetilde{\alpha_3} \in dcl_0(C_{2, \mathcal{S}} \beta_2 \beta^{\prime\prime\prime}_1 \alpha_1 \alpha^{\prime\prime\prime}_2 )$ and $\widetilde{\alpha_3} \in dcl_0(C_{2, \mathcal{S}} \beta_1 \beta^{\prime\prime\prime}_2 \alpha_2 \alpha^{\prime\prime\prime}_1 )$. Also by construction, we have $\alpha_1 \in dcl_0(C_{2, \mathcal{S}}  \widetilde{\alpha_3} \beta_2)$ and $\alpha^{\prime\prime\prime}_2 \in dcl_0(C_{2, \mathcal{S}}  \widetilde{\alpha_3} \beta^{\prime\prime\prime}_1)$, so  $\alpha_1 \alpha^{\prime\prime\prime}_2 \in dcl_0(C_{2, \mathcal{S}}  \widetilde{\alpha_3} \beta_2\beta^{\prime\prime\prime}_1)$. Similarly, we have 
				$\alpha_2 \alpha^{\prime\prime\prime}_1 \in dcl_0(C_{2, \mathcal{S}}  \widetilde{\alpha_3} \beta_1\beta^{\prime\prime\prime}_2)$, which concludes the proof of the claim.
			\end{proof}

			To keep similar notations, let $(\alpha_1, \alpha_2, \alpha_3, \beta_1, \beta_2, \beta_3)$ denote the tuple $(c_1, c_2, c_3, d_1, d_2, d_3)$ given by the claim. Now that we have the stronger interdefinability hypothesis, we can follow the proof of the usual group configuration theorem, to build the group $\Gamma$ as a group of germs of definable bijections. See \cite{bouscaren-config} or Definitions 3.1, 3.4, 3.7 in \cite[Version 3]{Wang_group_config} for the details. So, there exists some quantifier-free-$\mathcal{L}_{0}(C_2)$-type-definable group $\Gamma$ and elements $g_1, g_2, g_3 \in \Gamma(N_0)$ such that the group configurations drawn at the end of the statement of the theorem are equivalent over $M_1$.

			Our only remaining concerns are checking that we can replace $\Gamma$ by an $\mathcal{L}_{0}(M_1)$-definable group $\Gamma_0$ containing $\Gamma$, and that the group configuration built at the end of the proof can be made of elements of $N_1$. First note that the type-definable group $\Gamma$ is built as a group of germs of definable bijections acting on the stationary type $tp_0(\alpha_2 / M_1)$, said definable bijections coming from the interdefinability properties of the configuration given by Claim \ref{claim_wma_interdcl}. This implies that $\Gamma$ is  type-definable over $C_2$, in the sense of $T_0$.
			
			
			%
			%


			Let us now explain in more detail which elements of $N_1$ are used to build a group configuration for $\Gamma$.  Let $\beta^{\prime}_1, \beta^{\prime}_2, \beta^{\prime}_3 \in M_1$ be such that $\beta^{\prime}_1 \beta^{\prime}_2 \beta^{\prime}_3 \equiv_{C_2} \beta_1 \beta_2 \beta_3$. Then, the elements $g_1, g_2, g_3$ are all in the $\mathcal {L}_0$-quantifier-free-definable closure of $C_2 \beta^{\prime}_1 \beta^{\prime}_2 \beta^{\prime}_3 \alpha_2 \beta_1 \beta_2$. See Subsection 3.3 in \cite[Version 3]{Wang_group_config} for the construction.

			Finally, let us explain why we may assume the group to be definable rather than type-definable, while preserving all the required properties. First, any type-definable stable group is an intersection of definable groups (see for instance \cite[Proposition 3.4]{HruRK-MetaGp}). Then, by superstability, there are no infinite descending chains of definable subgroups where each one is of infinite index in the previous one (see for instance Proposition 1.5 and Theorem 5.18 in \cite{poizat-stable-groups}). By saturation of $M_{1, \mathcal{S}}$, the group $\Gamma$ is an intersection of quantifier-free-$\mathcal{L}_{0}(M_1)$-definable groups. Thus, there is some quantifier-free-$\mathcal{L}_{0}(M_1)$-definable group $\Gamma_0 \geq \Gamma$ such that $\Gamma$ is of bounded index inside $\Gamma_0$. Then, the generic elements $g_i$ in $\Gamma$ are also generic in $\Gamma_0$ over $M_{1, \mathcal{S}}$. \end{proof}

Combining this result with Theorem \ref{theo_how_to_find_group_homomorphisms}, we get the following

\begin{theo}\label{theo_group_with_kernel_blind_to_S}
	
	Let $T_0, T_1$ be theories, in languages $\mathcal{L}_0 \subseteq \mathcal{L}_1$ respectively. Assume that $T_0$ is superstable, has quantifier elimination and elimination of imaginaries in $\mathcal{L}_0$, in a collection of sorts $\mathcal{S}$. Assume that $T_1$ is $\mathrm{NIP}$. Let $N_0 \models T_0$, $N_1 \models T_1$, and $\iota : N_1 \rightarrow N_0$ be an $\mathcal{L}_0$-embedding. Assume that $dcl_0(\iota(N_1)) = \iota(N_1)$ and $acl_0(A) \cap \iota(N_1) = acl_1(A)$ for all $A \subseteq \iota(N_1)$.
	
	Let $V$ be a product of sorts of $\mathcal{S}$. Let $G$ be a definably amenable group, definable in $N_1$, with an $N_1$-definable map $\pi: G \rightarrow V$, whose fibers are blind to $\mathcal{S}$. Then, there exists a quantifier-free-$\mathcal{L}_0(N_1)$-definable group $H$ in the sorts of $\mathcal{S}$, and an $\mathcal{L}_1(N_1)$-definable group homomorphism $G^{00}_{N_1}(N_1) \rightarrow H(N_1)$ whose kernel is blind to the sorts in $\mathcal{S}$.
	
\end{theo}

\begin{proof}
	We may assume that $N_1$ and $N_0$ are very saturated. Let $M \preceq^+ N_1$ be a sufficiently saturated model over which everything is defined. Then, by Proposition \ref{prop_exists_generics}, there exists a type $p \in S_G(N_1)$ which is  strongly bi-f-generic  over $M$. Let $a \models p|_M$, $b \models p|_{Ma}$, $c \models p|_{Mab}$.


	Then, by Theorem \ref{theo_relative_stable_group_config}, there exists a quantifier-free-$\mathcal{L}_0(N_1)$-definable group $H$ and elements $\alpha, \beta$ in $H(N_1)$ that are independent and generic over $M$, such that the triples $(\alpha, \beta, \beta \cdot \alpha)$ and $(\pi(a), \pi(b), \pi(b \cdot a))$ are equivalent over $M$.

	We wish to conclude using Theorem \ref{theo_how_to_find_group_homomorphisms}. Let $\nu$ be the ideal on $G$ containing the sets blind to $\mathcal{S}$. Then, by construction, the type $tp(a /M\alpha)$ is in $\nu$. Also, by Proposition \ref{prop_the_blinders_are_a_nice_ideal}, the ideal $\nu$ is $\vee$-definable, $M$-invariant, invariant under left and right translations, and admits non-forking descent. Thus, we may apply Theorem \ref{theo_how_to_find_group_homomorphisms} in this setting, which yields an appropriate group homomorphism. Note that modding out by a finite normal subgroup is not an issue here, since $T_0$ eliminates imaginaries in the sorts of $\mathcal{S}$.
\end{proof}

\end{subsection}

\end{section}

\begin{section}{Algebraically bounded $\mathrm{NTP_2}$ Fields}

	Recall that a theory of (possibly enriched) fields $T$ is a theory of \emph{algebraically bounded fields} if for all $M \models T$, for all parameter sets $A \subseteq K(M)$, the algebraic closure $acl^M(A) \cap K$ coincides with the field-theoretic relative algebraic closure $K(A)^{alg} \cap K(M)$.
	
	

	\begin{defi}\label{defi_pure_im}
		Let $T$ be a theory of (possibly enriched) fields, and $M \models T$. Let $X$ be a set definable in $M^{eq}$. We say that $X$ is \emph{purely imaginary} if it is blind to $K$.

		We say that a field has a \emph{bounded Galois group} if, for all $n$, it has only finitely many Galois extensions of degree $n$ up to isomorphism. Equivalently, its absolute Galois group has only finitely many closed quotients of order $n$ for all $n$.

	\end{defi}

	With these definitions, the results of the previous section yield the following theorem:

	\begin{theo}\label{theo_general_morphism_with_imaginary_kernel}
		
		Let $T$ be an algebraically bounded theory of possibly enriched, possibly many-sorted, perfect fields. Assume that $T$ is $\mathrm{NIP}$, or that it is $\mathrm{NTP_2}$ and its models have bounded Galois groups. Let $G$ be a definably amenable group, admitting a definable map to some $K^n$, whose fibers are blind to $K$, in some sufficiently saturated model $M$ of $T$. Then, there exist an algebraic group $H$ over $M$ and an $M$-definable homomorphism $G^{00}_M \rightarrow H$, whose kernel is purely imaginary.
	\end{theo}
	
	\begin{proof}
		
		In the case where $T$ is $\mathrm{NIP}$, this is a consequence of Theorem \ref{theo_group_with_kernel_blind_to_S}, applied to $T_0 = \mathrm{ACF}$. So, let us assume that $T$ is $\mathrm{NTP_2}$ and its models have bounded Galois groups.

		Let $N$ be a very saturated elementary extension of $M$. We shall work inside $N$ and $N^{alg}$. First, using Proposition \ref{prop_how_to_build_group_configs_in_stable_reducts}, we construct a tuple $\overline{a b}$ of bi-f-generic elements of $G$, whose images in the field sort yield a group configuration, in the sense of $\mathrm{ACF}$, over $M$. Then, applying the usual group configuration theorem, there exists an algebraic group $H_1$ over $M^{alg}$ and a tuple of elements $\overline{\alpha \beta}$ in $H_1(N^{alg})$ which is equivalent over $M^{alg}$ (in the sense of $ACF$) with the field sort coordinates of $\overline{a b}$.

		Then, using boundedness of the Galois group of $M$, we may interpret, \emph{using parameters in $K(M)$}, the finite extensions of $M$ and $N$ generated by the parameters of definition of $H_1$ and the coordinates of $\overline{\alpha \beta}$ respectively. To fix notations, let $L_1 \geq M$ and $L_2 = N \cdot L_1$ denote these finite extensions. Such an expression of $L_2$ as a compositum exists by boundedness of the absolute Galois group of $K(M)$.
		
		Through the interpretations, we can see $H_1(L_2)$ as the group of $N$-points of an $\mathrm{ACF}$-definable group $H_2$, defined over $M$. Thus, we can interpret $\overline{\alpha \beta}$ as a tuple of elements of $H_2(N)$. In other words, we may assume that $H_1$ is over $M$, and that the points $\overline{\alpha \beta}$ are in $N$.

		Finally, we apply Theorem \ref{theo_how_to_find_group_homomorphisms}, to get some finite $H_0 \leq H_2(M)$, and an $M$-definable homomorphism $G^{00}_M \rightarrow Z_{H_2}(H_0) / H_0$ with purely imaginary kernel, where $Z_{H_2}(H_0) \leq H_2$ denotes the centralizer of $H_0$. Note that $Z_{H_2}(H_0)$ is quantifier-free definable, since $H_0$ is finite and $H_2$ is an algebraic group. Then, by elimination of imaginaries, the quotient $Z_{H_2}(H_0) / H_0$ is also quantifier-free $\mathcal{L}_{ring}(M)$-definable, so is $M$-definably isomorphic to an algebraic group over $M$. This concludes the proof.\end{proof}

	The aim is then to exploit the theorem above, to classify interpetable fields. Let us begin by proving some general results on intrepretable fields, and groups of affine transformations of such.

	\begin{lemma}\label{lemma_A=F}
		Let $F$ be an infinite field, definable in some arbitrary theory $T$. Let $B \leq F^{\times}$ be a type-definable subgroup of bounded index. Let $A \subseteq F$ be a $\vee$-definable subring of $F$ containing $B$. Then $A=F$.

	\end{lemma}

	\begin{proof}
		First, as $A$ contains $B$, we know that $F$ is covered by boundedly many multiplicative translates of $A$. By compactness, finitely many multiplicative translates of $A$ cover $F$. Thus, $F$ is a finitely generated $A$-module, and even a finite union of $A$-lines. 
		
		\begin{claim}
			The ring $A$ is a subfield of $F$. 
		\end{claim}
		\begin{proof}
			Let $a \in A$. We wish to show that $\frac{1}{a} \in A$. Multiplication by $\frac{1}{a}$ is an endomorphism of the finitely generated $A$-module $F$. Let $u : F \rightarrow F$ denote this endomorphism. Then, using the theorem of Cayley-Hamilton, see \cite[Chapter XIV, Section 3, Theorem 3.1]{Lan} or \cite{Stacks}[Part 1, Chapter 10, Lemma 10.16.2], there exists a monic polynomial $P \in A[X]$ such that $P(u) =0$. In particular, we have $P(u)(1) = 0$. So, some equality of the form $\frac{1}{a^n} + \frac{1}{a^{n-1}}c_1 + \cdots + c_n = 0$ holds, where the $c_i$ are in $A$. Multiplying by $a^{n-1}$, we deduce that $\frac{1}{a} \in A$, as required.
		\end{proof}
		
		Thus, $F$ is a finite extension of $A$, and even a finite union of $A$-lines. We know that both fields are infinite. Thus, by a well known fact on finite-dimensional vector spaces over infinite fields, $F$ is a single $A$-line, i.e. $F=A$.
	\end{proof}

	\begin{prop}\label{prop_facts_on_aff_1}
		Let $F$ be an infinite field, definable in some arbitrary theory $T$. Let $G$ denote the semi-direct product ${F_+} \rtimes {F^{\times}}$. 
		
		\begin{enumerate}
			\item Let $N$ be a normal subgroup of $G$, not contained in $F_+$. Then, $N$ contains $F_+$.
			\item Let $G_1 \leq G$ be a type-definable subgroup of bounded index. Then $G_1$ contains $F_+$, and $G_1$ has no nontrivial finite normal subgroups.

		\end{enumerate}

	\end{prop}
	
	\begin{proof}
		We begin with the following computation, left to the reader:
		
		\begin{claim}
			Let $(b, a)$, $(d, c)$ be elements of $G = {F_+} \rtimes {F^{\times}}$. Then, the following equality holds:  $(b, a) \cdot (d, c) \cdot (b, a)^{-1} = (b(1-c) + ad, c)$.
		\end{claim}

		Let us prove the first point. Let $(d, c) \in N$ be such that $c \neq 1$. Then, by the computation above, conjugating by elements $(b, 1)$ for $b \in F_+$, we get that all elements of the form $(e, c)$ belong to $N$, for $e \in F_+$. Multiplying by $(d,c)^{-1}$, we deduce that all elements of $F_+$ are in $N$.

		Let us now prove the second point. By definition, there exists some set $A$ such that $G_1$ contains $G^{00}_A$. By the third point of Fact \ref{fact_G^00}, the latter is normal in $G$. Since the subgroups $G^{00}_A$ and $G_1$ are of bounded index in $G$, their projections are of bounded index in $F^{\times}$, in particular, they are infinite subgroups. By the first point applied to $G^{00}_A$, this implies that $G^{00}_A$, and a fortiori $G_1$, contains $F_+$. Finally, let $N_0 \lhd G_1$ be a finite normal subgroup. By the first point, we know that $N_0$ is contained in $F_+$. Then, let $(d, 1) \in N_0$. Since $G_1$ projects onto an infinite subgroup of $F^{\times}$, and $N_0$ is finite, the claim above implies that $d=0$, hence $N_0$ is trivial, as required.
	\end{proof}

	\begin{prop}\label{prop_fields_whose_affine_group_embeds_in_an_algebraic_group_are_known}
		Let $\mathcal{L}$ be a language, and let $T$ be an $\mathcal{L}$-theory. Let $F$ and $K$ be infinite fields definable in some sufficiently saturated model $M$ of $T$. Assume that there is an algebraic group $H$ over $K(M)$ and an $M$-definable embedding $ \iota : {({F_+} \rtimes {F^{\times}})}^{00}_M \hookrightarrow H$. Then, $F$ admits an $M$-definable field embedding into some finite extension $L$ of $K$.

	\end{prop}
	
	Recall that any finite extension of a field is definable with parameters, in the ring language, inside some power of the field sort.

	\begin{proof} To simplify notations, let $G_1 = {({F_+} \rtimes {F^{\times}})}^{00}_M$, and $B = {F^{\times}}^{00}_M$. Let $A \subset M$ be a small parameter set over which everything is defined.
		\begin{claim}
			The group $G_1$ contains $F_+ \rtimes B$.
		\end{claim}
		
		\begin{proof}
			Since $G_1$ is of bounded index in $F_+ \rtimes F^{\times}$, Proposition \ref{prop_facts_on_aff_1} (2) implies that it contains $F_+$. Also, the bounded index subgroup $G_1$ projects onto a subgroup of bounded index in $F^{\times}$. Thus, the projection of $G_1$ contains $B$. This implies that $G_1$ contains $F_+ \rtimes B$. \end{proof}
		

		Now, we choose an embedding $K(M) \hookrightarrow {K(M)}^{alg}$. Let $\overline{U}$ denote the $K(M)^{alg}$-Zariski closure, inside $H$, of $\iota(F_+)$. Similarly, let $\overline{B}$ denote the $K(M)^{alg}$-Zariski closure of $\iota(B)$ inside $H$.

		\begin{claim}
			The set $\overline{U}$ is $K(M)^{}$-Zariski-closed. Similarly, the group $\overline{B}$ is $K(M)^{}$-Zariski-closed.
			
		\end{claim}
		\begin{proof}
			Let $\sigma \in Gal({K(M)}^{alg} / K(M)^{ins})$. Then, $\sigma$ fixes $\iota(F_+(M))$ pointwise, so $\sigma$ fixes the closure $\overline{U}({K(M)}^{alg})$ setwise. Since this holds for any $\sigma$, we have proved that $\overline{U}$ is $K(M)^{ins}$-Zariski-closed. Applying a suitable iterate of the Frobenius automorphism, we deduce that $\overline{U}$ is $K(M)^{}$-Zariski-closed. The same proof applies to $B$ and $\overline{B}$.
		\end{proof}
		
		Note that the claim implies that $\overline{U}$ is definable in $T$. 
		

		%
		%
	
	\begin{claim}
		The Krull dimension of $\overline{U}$ is equal to the maximal transcendence degree $d$, over $A$, of an element of $\iota(F_+)(M)$.
	\end{claim}
	\begin{proof}
		This is a well-known fact, which relies on compactness, which can be applied because $F_+$ corresponds to a clopen, in particular closed, set in the type space. Note that the result is false for $\vee$-definable sets.
	\end{proof}

	Then, since $B$ is of bounded index in $F^{\times}$, there is a bounded family $(a_i)_{i \in I}$ of elements of $F(M)$ such that $F = \bigcup_i B \cdot a_i$, where the dot denotes multiplication. Then, recall that, in the semidirect product $F_+ \rtimes B$, the action of $B$ on $F_+$, which coincides with multiplication in the field $F$, is by conjugation. So, for some $i$, the type-definable orbit ${\iota(a_i)}^{\iota(B)}$ contains an element of maximal transcendence degree among those in $\overline{U}$. Thus, the $ACF$-definable set $\iota(a_i)^{\overline{B}}$ is generic in the algebraic group $\overline{U}$.

	\emph{Note}: in general, the dimensions of $\overline{U}$ and $\overline{B}$ are not equal, because the embedding $F_+ \rtimes B \hookrightarrow H$ may be far from algebraic. For instance, in $DCF_0$, one might consider $(a, \lambda) \mapsto (a, \lambda, \frac{D\lambda}{\lambda})$.
	
	\begin{claim}\label{claim_V_is_definable}
		The subgroup $V$ of $\overline{U}$ generated by $\iota(a_i)^{\overline{B}}$ is definable and of finite index.
	\end{claim}
	
	\begin{proof}
		Since $\iota(a_i)^{\overline{B}}$ is generic, it suffices to prove definability. First, recall that the ACF-connected component, $\overline{U}^0$ of $\overline{U}$ exists and has finite index in $\overline{U}$ (see Fact \ref{fact_G^00}). Also, since some generic $p$ of $\overline{U}$ concentrates on $\iota(a_i)^{\overline{B}}$, we know by Remark \ref{rem_generics_in_stable_groups} that $\overline{U}^0$ is contained in $\iota(a_i)^{\overline{B}} \pm \iota(a_i)^{\overline{B}}$. Thus, the subgroup $V$ contains $\overline{U}^0$, and hence is a finite union of cosets of $\overline{U}^0$. This shows definability.
	\end{proof}
	
	
	Recall that $\iota(B)$ acts on $\iota(F_+)$ by conjugation. Then, one can check that the closure $\overline{B}$ acts on the closure $\overline{U}$ by conjugation as well. Hence, the algebraic group $\overline{B}$ acts on the orbit $\iota(a_i)^{\overline{B}}$ via group homomorphisms, so it acts on the group generated by this orbit, which is $V$. Then, let $R$ be the ring of endomorphisms of the group $V$ generated by $\overline{B}$. It is a priori $\vee$-definable in $M^{alg}$.

	\begin{claim}\label{claim_R_commutative}
		The ring $R$ is commutative, and its action on $V$ is $\overline{B}$-equivariant.
	\end{claim}
	%
	\begin{proof}
		First note that, as $B$ is abelian, so is $\iota(B)$. Thus, the Zariski closure $\overline{B} \leq H$ is also an abelian group. Since $R$ is generated by $\overline{B}$, it suffices to prove that any element of $R$ commutes with all elements of $\overline{B}$. So, let $u = b_1 \pm \cdots \pm b_n \in R$, with $n < \omega$, and the $b_i$ in $\overline{B}$, and let $b \in \overline{B}$. For any $v \in V$, we have $b(u(v)) = b(b_1(v)) \pm \cdots \pm b(b_n(v)) = (b_1 \cdot b) (v) \pm \cdots \pm (b_n \cdot b) (v) = u(b(v))$. This proves the result. 
	\end{proof}

	\begin{claim}
		The ring $R$ is definable (possibly with parameters in $M^{alg}$), and so is its action on $V$.
		
	\end{claim}
	
	\begin{proof}
		Let us show definability of the ring. By Claim \ref{claim_V_is_definable} above, and compactness, there exists some $N$ such that $V = {(\iota(a_i)^{\overline{B}})}^{\pm N}$, i.e. any element of $V$ can be written as a signed sum of at most $N$ elements of  $\iota(a_i)^{\overline{B}}$. Also note that, by Claim \ref{claim_R_commutative}, any element $f$ of the ring $R$ is uniquely determined by the image $f(\iota(a_i))$. Thus, we may encode the ring $R$ as a collection of $N$-tuples of elements of $\iota(a_i)^{\overline{B}}$, along with $N$-tuples of signs. 
	\end{proof}

	Thus, as we know that $ACF$-definable functions are given piecewise by rational fractions, there is a finite extension $L_1$ of $K(M)$ such that, for all finite extensions $L_2$ of $L_1$, the set $R(L_2)$ is a commutative subring of $R$. We may also assume that the image of $B(K(M))$ under the $ACF$-definable group embedding $\overline{B} \hookrightarrow R^{\times}$ is contained in $R(L_1)$.
	
	\begin{claim}
		There exists a finitely generated proper ideal $I$ of $R$ such that $R/I$ is definably isomorphic, in $ACF$, to the field $K(M)^{alg}$.
	\end{claim}
	
	\begin{proof}
		The idea is to mod out by zero divisors in the successive quotients, until one gets a domain $R/I$. First note that the group of units $R^{\times}$ is infinite, for it contains a copy of $B$. Hence, any nonzero definable ideal is infinite, thus has Morley rank at least $1$. Then, by additivity of the rank, modding out $R$ by a nontrivial definable ideal $I$ decreases the rank. Thus, after finitely many steps, the quotient ring $R/I$ has no zero divisors, i.e. is a domain.

		This quotient is then an $ACF$-definable domain, whose field of fractions is definably isomorphic to the field $M^{alg}$, by \cite[Theorem 4.15]{poizat-stable-groups}. As $K(M)^{alg}$ has no infinite definable subrings, we have shown that $R/I$ is definably isomorphic to $K(M)^{alg}$, as required.
	\end{proof}
	
	Let $R_1$ denote the $ACF$-definable ring $R/I$. Then, let $L_2$ be a finite extension of $L_1$ containing the parameters defining $R_1$ and the quotient map $R \rightarrow R_1$, and such that $R_1(L_2)$ is a subring of $R_1$. We shall now work inside the field $L_2$, seen as a definable set in $M$, in several variables. We now consider the following definable subset of $F$:  $A := \lbrace \lambda \in F \, | \, \exists u \in R(L_2) \, u|_{\iota(F_+) \cap V} = \lambda|_{\iota(F_+) \cap V} \rbrace$, where any $\lambda \in F$ acts on $F_+$, thus on $\iota(F_+)$, by multiplication.

	\begin{claim}
		The set $A$ is a subring of $F$.
	\end{claim}
	
	\begin{proof}
		The elements $0$ and $1$ of $R$ are indeed in the subring $R(L_2)$, and they act as $0$ and the identity respectively. The fact that $A$ is an additive subgroup of $F$ follows from the fact that $R(L_2)$ is an additive subgroup of $R$. Now, let us prove that $A$ is closed under multiplication. Let $\lambda_1, \lambda_2 \in F$, and $u_1, u_2 \in R(L_2)$ such that $\lambda_i$ and $u_i$ agree on $\iota(F_+) \cap V$, for $i=1,2$. Then, as $\lambda_2$ fixes $\iota(F_+)$ setwise, and $u_2$ fixes $V$ setwise, and since they agree on the intersection, we know that $u_2$ fixes $\iota(F_+) \cap V$ setwise. Then, we compute that $u_1 \circ u_2|_{\iota(F_+) \cap V} = \lambda_1 \circ u_2|_{\iota(F_+) \cap V} = \lambda_1 \cdot \lambda_2|_{\iota(F_+) \cap V}$. This show that $\lambda_1 \cdot \lambda_2 \in A$, as required.
	\end{proof}

	\begin{claim}
		The subring $A$ is equal to $F$. 
	\end{claim}
	
	\begin{proof}
		The subring $A$ is definable. Also, by construction, the subring $A$ contains $B$, which is of bounded index. Thus, by Lemma \ref{lemma_A=F}, we have $A=F$. 
	\end{proof}

	Hence, there is a ring embedding $F=A \hookrightarrow R(L_2)$, which is definable in $M$. Finally, composing with the quotient map $R \rightarrow R_1$, we get a definable ring morphism $F \rightarrow R_1(L_2)$, which is an embedding since $F$ is a field. Since $R_1$ is in fact $K(M)^{alg}$ (i.e. the tautological definable ring), we have found an $M$-definable embedding $F \hookrightarrow L_2$, with $L_2$ interpreted in $M$ in the usual way.
\end{proof}

\begin{lemma}\label{lemma_actually_embed_affine_group}
Let $T$ be a theory of possibly enriched, possibly many-sorted, fields. Let $M$ be a sufficiently saturated model of $T$, $F$ an infinite field interpretable in $M$, $H$ an algebraic group over $M$ and $f : {(F_+ \rtimes F^{\times})}^{00}_M \rightarrow H$ an $M$-definable group homomorphism with purely imaginary kernel. Assume that $F$ is not purely imaginary. Then, the morphism $f$ is injective.

\end{lemma}

\begin{proof}
Let $A$ denote the purely imaginary kernel of $f$. Also, let $G$ denote the group ${(F_+ \rtimes F^{\times})}^{00}_M $. We wish to show that $A$ is trivial. We know by Proposition \ref{prop_facts_on_aff_1} that $G$ contains $F_+$.

\begin{claim}
The subgroup $A$ is contained in $G \cap F_+ = F_+$.
\end{claim}

\begin{proof}
By contradiction, assume that $A$ is not contained in $F_+$. Let $(b, a) \in A \setminus F_+$. Then, for all $(d, 1) \in F_+$, we have $(b, a)^{(d,1)} = (b+ (1-a) d, a) \in A$. Since $a \neq 1$, we just found a definable injection of sets  $F_+ \hookrightarrow A$.  Since $A$ is purely imaginary, but $F$ is not, this is a contradiction.
\end{proof}


Thus, the restriction of $f$ to $G \cap F^{\times}$ is a definable embedding $G \cap F^{\times} \hookrightarrow H$. Then, by compactness, there exists a definable set $Y$ contained in $F$ such that $ G \cap F^{\times} \subseteq Y$ and a definable injection of sets $Y \hookrightarrow H$. Thus, $Y$ definably embeds into some power of the field sort. Also, by compactness, the set $F$ is a finite union of multiplicative translates of $Y$. So, the definable set $F$ also embeds into some power of the field sort. Hence, the type-definable set $G$ definably embeds into some power of the field sort.

Since $A$ is purely imaginary and embeds into some power of the field sort, it is thus finite. By Proposition \ref{prop_facts_on_aff_1}, the group $G$ admits no nontrivial finite normal subgroups. Thus $A$ is trivial, and we are done.
\end{proof}


\begin{coro}\label{coro_fields_are_either_pur_im_or_algebraic}
Let $T$ be an algebraically bounded theory of possibly enriched, possibly many-sorted, perfect fields. Assume that $T$ is $\mathrm{NIP}$, or that it is $\mathrm{NTP_2}$ and its models have bounded Galois groups. Let $F$ be an infinite definable field, admitting a definable map to some $K^n$, whose fibers are blind to $K$, in some sufficiently saturated model $M$ of $T$. Then either $F$ is purely imaginary, or $F$ admits a definable embedding into some finite extension of $K(M)$.

\end{coro}

\begin{proof}
Assume that $F$ is not purely imaginary. We apply Theorem \ref{theo_general_morphism_with_imaginary_kernel}, to the solvable, thus definably amenable, group $F_+ \rtimes F^{\times}$. So, we find an algebraic group $H_1$ over $M$, and a definable homomorphism $g : {(F_+ \rtimes F^{\times})}^{00}_M \rightarrow H_1$, with purely imaginary kernel. By Lemma \ref{lemma_actually_embed_affine_group}, the morphism $g$ is an embedding. The result then follows from Proposition \ref{prop_fields_whose_affine_group_embeds_in_an_algebraic_group_are_known}.
\end{proof}

%
%
%
%
%
%
%
%

\begin{ex}\label{ex_classes_of_examples}
Theorem \ref{theo_general_morphism_with_imaginary_kernel} and Corollary \ref{coro_fields_are_either_pur_im_or_algebraic} above cover the following examples: NIP or $\mathrm{NTP}_2$ bounded henselian valued fields of characteristic $0$, for which we shall prove the more specific Theorem \ref{theo_fields_hen_0}, bounded pseudo real-closed fields (see \cite[Section 4]{Stabilizers-NTP2}), bounded pseudo $p$-adically closed fields (see \cite{Mon-NTP2}), and bounded PAC fields. Also note that all the examples mentioned in the introduction are covered by these results.
\end{ex}

\end{section}

\begin{section}{Examples}
	
	\begin{subsection}{Algebraically bounded differential fields}

		Here, we introduce a tame class of theories of differential fields, inside which our general tools can be used to study intepretable groups and fields. The examples are essentially the same as those studied in \cite{ppp_def_groups_codf}, see Remark \ref{rem_what_about_examples} below. Recall that ``purely imaginary'' is synonym of ``blind to the field sort''. If $A$ is a subset of a differential field, we let $\lbrace A \rbrace_{.^{-1}, \partial}$ denote the differential subfield generated by $A$.

		\begin{defi}\label{defi_alg_bounded_diff_fields}
			Let $T$ be a theory of differential fields, in a language extending the language of differential rings, possibly many-sorted. We say that $T$ is a theory of \emph{algebraically bounded differential fields} if, for all $M \models T$, for all $A \subseteq K(M)$, the model-theoretic algebraic closure of $A$ inside $M$ coincides with the field-theoretic relative algebraic closure of $\lbrace A \rbrace_{.^{-1}, \partial}$.
			
		\end{defi}

		\begin{theo}\label{theo_groups_in_algebraically_bounded_differential_fields}
			
			Let $T$ be an NIP theory of algebraically bounded differential fields. Let $M$ be a sufficiently saturated model of $T$. Let $n < \omega$.
			\begin{enumerate}
				\item Let $G$ be an $M$-definable group, admitting an $M$-definable map $G \rightarrow K^n$ whose fibers are blind to $K$. Assume that $G$ is definably amenable. Then, there exists an $M$-definable group homomorphism $G^{00} \rightarrow H$ with kernel blind to $K$, where $H$ is an algebraic group over $K(M)$.
				
				\item Let $F$ be an infinite $M$-definable field, admitting an $M$-definable map $F \rightarrow K^n$ whose fibers are blind to $K$.  Then, either $F$ is blind to $K$, or $F$ admits a definable embedding into some finite extension of $K(M)$.
				
			\end{enumerate}
			
		\end{theo}
		\begin{proof}
			Let us prove the first point. We wish to apply Theorem \ref{theo_group_with_kernel_blind_to_S} with $T_0 = \mathrm{DCF_{0}}$ and $T_1 = T$, and the collection $\mathcal{S}$ containing only the field sort itself. The definition of algebraically bounded differential fields implies that the technical hypotheses of Theorem \ref{theo_group_with_kernel_blind_to_S} hold. A priori, we get a group $H$ which is $\mathcal{L}_0(M)$-quantifier-free definable, where $\mathcal{L}_0 = \mathcal{L}_{rings} \cup \lbrace \partial \rbrace$, and an $M$-definable morphism $G^{00}(M) \rightarrow H(M)$ with purely imaginary kernel. Then, by \cite[Corollary 4.2]{Pillay1997SomeFQ}, the group $H$ actually embeds $M$-definably  into an algebraic group over $M$, which proves the first point.

			Then, let us prove the second point. Appyling the first point to the solvable group $F_+ \rtimes F^{\times}$, we get an algebraic group $H$ over $M$, and an $M$-definable morphism $g: {(F_+ \rtimes F^{\times})}^{00} \rightarrow H$ with kernel blind to $K$. If $F$ is not blind to $K$, then by Lemma \ref{lemma_actually_embed_affine_group}, the morphism $g$ is injective. We can then conclude using Proposition \ref{prop_fields_whose_affine_group_embeds_in_an_algebraic_group_are_known}. \end{proof}

\end{subsection}

\begin{subsection}{Differentially closed valued fields}

\begin{subsubsection}{Some background on differentially closed valued fields}

	If $(K, v)$ is a valued field, there are several imaginary sorts, called the \emph{geometric sorts}, that are of interest. We let $\mathcal{O}$ denote the valuation ring, $\mathfrak{m}$ its maximal ideal, and $k$ the residue field $\mathcal{O} / \mathfrak{m}$. For each $n$, the sort $S_n$ is the quotient $GL_n(K) / GL_n(\mathcal{O})$, coding the $\mathcal{O}$-lattices of rank $n$, $\Lambda \leq K^n$. Also, for each $n$, the sort $T_n$ is the set of pairs $(\Lambda, u)$, where $\Lambda \in S_n$ is a lattice, and $u \in \Lambda / \mathfrak{m} \Lambda$ is an element of the $n$-dimensional $k$-vector space $\Lambda / \mathfrak{m} \Lambda$. So, we have a definable surjection $T_n \rightarrow S_n$, whose fibers are $k$-vector spaces of dimension $n$. Together with the sort $K$, and the appropriate quotient maps from powers of $K$ to $T_n$ and $S_n$, these define the \emph{geometric language}.


	Let us first recall a few facts about algebraically closed valued fields.
	
	\begin{fact}\label{fact_ACVF}
		\begin{enumerate}
			
			\item The theory of valued fields has a model completion $\mathrm{ACVF}$, the theory of non-trivially valued algebraically closed fields. The theory $\mathrm{ACVF}$ eliminates quantifiers in the three-sorted language (with sorts for $K$, $k$, $\Gamma$), and its completions are given by the characteristics of the valued field and the residue field. It is an algebraically bounded theory of fields.

			\item The residue field is a stably embedded pure algebraically closed field. The value group is a stably embedded pure divisible ordered abelian group.

			\item The theory $\mathrm{ACVF}$ has elimination of imaginaries in the geometric sorts \cite[Theorem 1.0.1]{HasHruMac-ACVF}, and is NIP.

		\end{enumerate}

	\end{fact}

	The following theorem classifies the interpretable fields in ACVF:
	
	\begin{theo}[See Theorem 6.23 in \cite{HruRK-MetaGp}]\label{theo_fields_ACVF}
		Let $F$ be an infinite field interpretable in ACVF. Then $F$ is definably isomorphic to the valued field or the residue field.
	\end{theo}

	Recall that the \emph{separant} $s_P$ of a multivariate polynomial $P$ is the polynomial $\frac{\partial P}{\partial X_n}$, where $n$ is the smallest index such that $P$ is in the variables $X_0,\cdots, X_n$. We define separants for differential polynomials similarly.

	\begin{fact}\label{fact_DCVF}
		\begin{enumerate}
			
			\item The theory of differential valued fields of equicharacteristic $0$ has a model completion $\mathrm{DCVF}$. The theory $\mathrm{DCVF}$ eliminates field sort quantifiers in the geometric language, and is complete.
			
			\item An axiomatization of $\mathrm{DCVF}$, inside the finitely axiomatizable class of differential valued fields, is as follows : for every differential polynomial $P(x) \in K \lbrace x \rbrace$ with $|x| = 1$ and $ord_x(P) = m \geq 1$, for field sort variables $y = (y_0,..., y_m)$, the following formula holds

			$\forall \gamma (\exists y (P^*(y) = 0 \wedge s_{P^*}(y) \neq 0)) \rightarrow \exists x (P(x) = 0 \wedge s_P(x) \neq 0 \wedge \bigwedge_i v(\partial^i x - y_i) \geq \gamma).$

			Here, $P^*$ denotes the ordinary polynomial associated to $P$, and $s_P$ denotes the separant of the differential polynomial $P$.

			\item The residue field is a stably embedded pure algebraically closed field. The value group is a stably embedded pure divisible ordered abelian group.

			\item The theory $\mathrm{DCVF}$ has elimination of imaginaries in the geometric sorts, and is NIP.

		\end{enumerate}
	\end{fact}
	
	Proofs for points 1 and 2 of the fact above can be found in \cite[Theorems 2.3.4 and 2.4.2]{CKPoi-DCVF}. For the elimination of imaginaries part of the fourth point, see  \cite[Theorem 3.3.2]{CKPoi-DCVF} or \cite[Theorems 8.7 and 9.3]{Rid-VDF}. 
	

	\begin{rem}\label{rem_stable_sorts_DCVF}
		
		Let $\mathcal{L}^{\mathcal{G}}$ be a Morleyisation of the geometric language for $\mathrm{ACVF}$, so that quantifier elimination holds. Then, let $\mathcal{L}^{\mathcal{G}}_{\partial} = \mathcal{L}^{\mathcal{G}} \cup \lbrace \partial \rbrace $ be the extended language, for $\mathrm{DCVF}$, where $\partial$ denotes the derivation on the field sort. Then, by Fact \ref{fact_DCVF}, the theory DCVF eliminates imaginaries in the language $\mathcal{L}^{\mathcal{G}}_{\partial}$.

	\end{rem}

	If $a$ is a finite tuple in a differential field, we let $\nabla a$ denote the infinite tuple of all the derivatives of the coordinates of $a$. Similarly, if $d < \omega$, we let $\nabla_d a$ denote the finite tuple of the derivatives of order $\leq d$ of the coordinates of $a$. We can also do this for tuples of variables, to define terms in the language of differential fields.

	\begin{prop}\label{prop_dcl_in_field_dcvf}
		
		Let $A \subseteq K(M)$, where $M \models DCVF$. Then, we have the following equality : $K(dcl(A)) = K(dcl^{ACVF}(\lbrace A \rbrace_{.^{-1}, \partial}))$. Similarly, $K(acl(A)) = K(\lbrace A \rbrace_{.^{-1}, \partial})^{alg}$.

	\end{prop}
	
	\begin{proof}
		Let us first prove the second point. Let $c \in K(acl(A))$. Let $\phi(x, a)$ be an $A$-formula with only finitely many realizations, where $|x|=1$, satisfied by $c$. By quantifier elimination (see for instance \cite[Corollary 2.4.7]{CKPoi-DCVF}), it is enough to deal with the case where $\phi(x,y)$ is of the form $\psi(\nabla x, \nabla y)$, and $\psi$ is an ACVF-formula without parameters, of the form $\bigwedge_i P_i(\overline{x}, \overline{y}) = 0 \wedge \theta(\overline{x}, \overline{y})$, where the $P_i$ are polynomials with integer coefficients, and $\theta$ defines an open set (in ACVF). In fact, using differential algebra (see \cite{marker-differential-fields}, Lemma 1.8), if $P$ is a minimal differential polynomial for $c$ over the field $\lbrace A \rbrace_{.^{-1}, \partial}$, we have $P_i(\overline{x}, \nabla a)$ belongs to the differential prime ideal defined by $P$ for all $i$, so the formula $\phi(x, a)$ is implied by the formula $P(x) = 0 \wedge s_P(x) \neq 0 \wedge \theta(\nabla x, \nabla a)$. Thus, the latter only has finitely many realizations, and is satisfied by $c$. Let $P^*$ be the polynomial such that $P^*(\nabla x) = P$. Then, using the axioms of DCVF, we can check that the formula $P^*(\overline{x}) = 0 \wedge  s_{P^*}(\overline{x}) \neq 0 \wedge \theta(\overline{x}, \nabla a)$ only has finitely many solutions in ACVF. Otherwise, consider an arbitrarily large finite number of solutions, separate them using balls, and find solutions which are prolongations of elements inside each of these balls.

		Now, let us prove the first point. Remember that we are in characteristic $0$. Note that the field $B=K(dcl^{ACVF}(\lbrace A \rbrace_{.^{-1}, \partial}))$ is Henselian, because it is $dcl$-closed in ACVF. So, its valuation extends uniquely, up to isomorphism, to its field-theoretic algebraic closure. Also, by standard computations in differential algebra, the derivation extends uniquely to the field-theoretic algebraic closure. Thus, the field-theoretic automorphisms of $B^{alg}$ fixing $B$ pointwise are actually automorphisms of the differential valued field structure. So, any element of $K(dcl^{DCVF}(B))$ is actually in $B$.
	\end{proof}

	The following definable types are useful for some manipulations.
	
	\begin{defi}
		Let $\Lambda \subseteq K^n$ be a lattice (so $\Lambda$ is coded by an element of $S_n$). We let $\eta^{ACVF}_{B\Lambda}$ denote the $ACVF$-generic type of $\mathcal{O}$-bases of $\Lambda$. It is a definable type. Then, by quantifier elimination for DCVF in the $3$-sorted language, we let $\eta^{c}_{B\Lambda}$ denote the type generated  by $\eta^{ACVF}_{B\Lambda}(x)$ and the formula $\bigwedge_i (\partial x_i = 0)$. This determines a complete definable type in the $3$-sorted language, which corresponds to $ACVF$-generic bases of $\Lambda$, whose elements have constant coordinates.

		Then, this definable type extends uniquely to a definable type in the extended language  $\mathcal{L}^{\mathcal{G}}_{\partial}$, since we are only adding imaginaries. We shall denote it $\eta^{c}_{B\Lambda}$ as well. 
	\end{defi}
	
	\begin{rem}
		The definable type $\eta^{c}_{B\Lambda}$ ``commutes with itself'' for the tensor product operation. In fact, it commutes with all invariant types.
	\end{rem}


\end{subsubsection}

\begin{subsubsection}{Relative quantifier elimination in the geometric sorts}
	
	We now wish to prove the following relative quantifier elimination result : 
	
	\begin{theo}\label{theo_relative_QE_DCVF}
		The theory $DCVF$ admits quantifier elimination in the language $\mathcal{L}^{\mathcal{G}}_{\partial}$.

	\end{theo}

	We shall use standard back-and-forth arguments.

	\begin{prop}\label{prop_embedding_DCVF_geom}
		
		Let $M$ be an $\aleph_1$-saturated model of $DCVF$, and $N \succeq^+ M$ be a very saturated model of $DCVF$. Let $A \leq M$ be a countable substructure, and $f: A \rightarrow N$ an $\mathcal{L}^{\mathcal{G}}_{\partial}$-embedding. Then, $f$ can be extended into an $\mathcal{L}^{\mathcal{G}}_{\partial}$-embedding $M \rightarrow N $.

	\end{prop}

	To prove this result, we wish to use the elimination of field quantifiers, and prove that one can lift elements of $S_n(dom(f))$ and $T_n(dom(f))$ to the field sort in such a way that the embedding can be extended to those lifts.

	\begin{prop}\label{prop_qf_defi_for_generic_constant_bases}
		Let $\phi(x,y)$ be a quantifier-free formula in $\mathcal{L}^{\mathcal{G}}_{\partial}$, and $\Lambda \in S^n$ be a lattice. Then, the definition $d_{\eta^{c}_{B\Lambda}} x \, \phi(x, y)$ is also quantifier-free in $\mathcal{L}^{\mathcal{G}}_{\partial}$. More precisely, it is an $\mathcal{L}^{\mathcal{G}}_{\partial}$-quantifier-free formula in $(y, \Lambda)$.

	\end{prop}
	
	\begin{proof}
		Let $p=\eta^{c}_{B\Lambda}$, and $p_0 = p|_{ACVF} = \eta^{ACVF}_{B\Lambda}$. By construction of the language $\mathcal{L}^{\mathcal{G}}_{\partial}$, and the Leibniz rule, there exists a formula $\psi(\overline{x}, \overline{y})$ which is quantifier-free in $\mathcal{L}^{\mathcal{G}}_{}$, such that $DCVF \models \phi(x,y) \leftrightarrow \psi(\nabla x, \nabla y)$. Then, we have the following: 
		
		$$DCVF \models [\phi(x,y) \wedge \bigwedge_i (\partial x_i = 0)] \leftrightarrow [\psi(x, 0, \cdots, 0, \nabla y)  \wedge \bigwedge_i (\partial x_i = 0)].$$
		
		So, if $c$ is a $y$-tuple, we have $\phi(x, c) \in p(x)$ if and only if $\phi(x,c) \wedge \bigwedge_i (\partial x_i = 0) \in p(x)$, if and only if $\psi(x, 0, \cdots, 0, \nabla c)  \wedge \bigwedge_i (\partial x_i = 0) \in p(x)$, if and only if the formula $\psi(x, 0, \cdots, 0, \nabla c)$ is in $p_0(x)$. So, the $p$-definition of $\phi(x,y)$ is equivalent to the formula $d_{p_0} x \, \psi(x, 0, \cdots, 0 , \nabla y)$. By assumption, $ACVF$ has quantifier elimination in $\mathcal{L}^{\mathcal{G}}$, so  the definitions of $p_0$ are quantifier-free in $\mathcal{L}^{\mathcal{G}}$. Thus, the formula $\chi(y) = d_{p_0} x \, \psi(x, 0, \cdots, 0 , \nabla y)$ is quantifier-free in  $\mathcal{L}^{\mathcal{G}}_{\partial}$.

		Now, let us prove that the definition is quantifier-free in $(y, \Lambda)$. Let $\theta(x, \overline{y})$ be the formula $\psi(x, 0, \cdots, 0 , \overline{y})$. Recall that $\Lambda$ and $\mathcal{O}^n$ are, in $ACVF$, generically stable connected definable sub-$\mathcal{O}$-modules of $K^n$, and that isomorphisms of $\mathcal{O}$-modules are parametrized by $\mathcal{O}$-bases of $\Lambda$. If $b$ is such a basis, let $f_b : \mathcal{O}^n \simeq \Lambda$ be the $b$-definable isomorphism given by $b$. Now, it is not hard to check that the $\theta(x, \overline{y})$-definition of the generic type $\eta_{B\Lambda}^{ACVF}$ is equivalent, in $ACVF$, to the following formula : $ (\forall f_b : \mathcal{O}^n \simeq \Lambda) \,  [d_{\eta^{ACVF}_{B\mathcal{O}^n}} \, x \, \theta(f_b(x), \overline{y})]$. By quantifier elimination for $ACVF$ in $\mathcal{L}^{\mathcal{G}}_{}$, this formula is equivalent to an $\mathcal{L}^{\mathcal{G}}_{}$-quantifier-free formula $\rho(\Lambda, \overline{y})$. Then, the formula $\rho(\Lambda, \nabla y)$ is quantifier-free in $\mathcal{L}^{\mathcal{G}}_{\partial}$, and is a $\phi(x,y)$-definition of  $\eta^{c}_{B\Lambda}$.
	\end{proof}

	\begin{prop}\label{prop_qf_defi_for_generics_of_open_balls}
		Let $\phi(x,y)$ be a quantifier-free formula in $\mathcal{L}^{\mathcal{G}}_{\partial}$, and $b$ be (the code for) an open ball. Let $\eta^c_b$ be the definable type generated by $\eta^{ACVF}_b(x) \cup \lbrace \partial x = 0 \rbrace$. Then, the definition $d_{\eta^{c}_{b}} x \, \phi(x, y)$ is an $\mathcal{L}^{\mathcal{G}}_{\partial}$-quantifier-free formula in $(y, b)$.

	\end{prop}
	
	\begin{proof}
		The proof is very similar to that of Proposition  \ref{prop_qf_defi_for_generic_constant_bases}. The only notable difference lies in the parametrization used to get a quantifier-free formula in $(y, b)$ : we parametrize definable bijections between open balls $b$ and the open ball $\mathfrak{m}$ via affine transformations $f_{c, d} : x \in \mathfrak{m} \mapsto dx + c \in B(c, v(d))$.
	\end{proof}

	\begin{lemma}\label{lemma_WMA_lattices_are_lifted}
		In the context of Proposition \ref{prop_embedding_DCVF_geom}, we may extend the embedding $f$ such that all elements of $S_n(dom(f))$ have bases in $K(dom(f))$, while keeping $dom(f)$ countable.
	\end{lemma}
	
	\begin{proof}
		First assume that we know how to lift one lattice. Then, we can construct an $\omega$-chain of embeddings $f_k : A_k \rightarrow N$, such that all the elements of the $S_n(dom(f_k))$ are lifted in $K(dom(f_{k+1}))$, and each $A_k$ is still countable.  
		
		So, let us prove that we can indeed lift one lattice. Let $\Lambda \in S_n(dom(f))$ be a lattice. We intend to lift $\Lambda$ via a realization $a$ of $\eta^c_{B\Lambda}|_{dom (f)}$, then lift $f(\Lambda)$ via a realization $b$ of $\eta^c_{Bf(\Lambda)}|_{im (f)}$, and then to use Proposition \ref{prop_qf_defi_for_generic_constant_bases} to prove that sending $a$ to $b$ defines an $\mathcal{L}^{\mathcal{G}}_{\partial}$-embedding.
		
		Since $M$ and $N$ are $\aleph_1$-saturated, finding $a$ and $b$ is straightforward. Now, let us check that extending $f$ by sending $a$ to $b$ defines an embedding. So, let $\phi(x, y)$ be an $\mathcal{L}^{\mathcal{G}}_{\partial}$-quantifier-free formula, and $c \in dom(f)$ a finite tuple. By  Proposition \ref{prop_qf_defi_for_generic_constant_bases}, let $\psi(y, u)$ be an $\mathcal{L}^{\mathcal{G}}_{\partial}$-quantifier-free formula such that $\psi(y, \Omega)$ is a $\phi(x,y)$-definition of $\eta_{B\Omega}^{c}$ for all lattices $\Omega \in S_n$. In particular, $\psi(y, \Lambda)$ is, in $M$, a $\phi(x,y)$-definition of $\eta_{B\Lambda}^{c}$, and $\psi(y, f(\Lambda))$ is, in $N$, a $\phi(x,y)$-definition of $\eta_{Bf(\Lambda)}^{c}$. So, by choice of $a$ and $b$, we have $M \models \phi(a, c)$ if and only if $M \models \psi(c, \Lambda)$, if and only if $N \models \psi(f(c), f(\Lambda))$, if and only if $N \models \phi(b, f(c))$. \end{proof}

	\begin{lemma}\label{lemma_WMA_k_is_lifted}
		With the additional assumption of Lemma \ref{lemma_WMA_lattices_are_lifted}, we may extend the embedding $f$ such that all elements of $T_n(dom(f))$ are lifted in $K(dom(f))$, while keeping $dom(f)$ countable.
		
	\end{lemma}
	
	\begin{proof}
		Assume that all the elements of $S_n(dom(f))$ have bases in $K(dom(f))$. In particular, for all $\Lambda \in S_n(dom(f))$, there are $dom(f)$-quantifier-free definable bijections $\Lambda / \mathfrak{m} \Lambda \simeq k^n$.  We now use Proposition \ref{prop_qf_defi_for_generics_of_open_balls}, to lift elements of the residue field using generic constants. Note that, for definable types, having quantifier-free definitions goes through taking tensor products.  As in Lemma \ref{lemma_WMA_lattices_are_lifted}, the construction provides us with embeddings in the language $\mathcal{L}^{\mathcal{G}}_{\partial}$.
	\end{proof}
	
	\begin{proof}[Proof of Proposition \ref{prop_embedding_DCVF_geom}]
		
		Applying Lemmas \ref{lemma_WMA_lattices_are_lifted} and \ref{lemma_WMA_k_is_lifted} successively in an $\omega$-chain, we may assume that both properties of Lemmas \ref{lemma_WMA_lattices_are_lifted} and \ref{lemma_WMA_k_is_lifted} hold simultaneously, i.e., that the domain of $f$ is generated by its field sort.
		
		Finally, given an $\mathcal{L}^{\mathcal{G}}_{\partial}$-embedding $f : A \rightarrow N$, where $A$ is generated by its field sort, we use quantifier elimination in the $3$-sorted language $\mathcal{L}_{k, \Gamma}$, to extend $f$ to an $\mathcal{L}_{k, \Gamma}$-\emph{elementary} embedding $g$ defined on $K(M)$. Checking that such a map $g$ extends (uniquely) to an $\mathcal{L}^{\mathcal{G}}_{\partial}$-embedding defined on $M$, which extends $f$, is left to the reader.
	\end{proof}

	
	\begin{coro}\label{coro_closures_in_dcvf}
		
		\begin{enumerate}
			\item Let $M \models DCVF$. Let $A$ be a substructure of $M$ in the language $\mathcal{L}^{\mathcal{G}}_{\partial}$. Then, for any geometric sort $S$, we have $S(dcl(A)) = S(dcl^{ACVF}(A))$ and $S(acl(A)) = S(acl^{ACVF}(A))$.
			
			\item The theory $DCVF$ is a theory of algebraically bounded differential fields.
		\end{enumerate}
	\end{coro}

	\begin{proof}
		For the field sort, we use Proposition \ref{prop_dcl_in_field_dcvf}, which also implies the second point. For the other sorts, we use the relative quantifier elimination result of Theorem \ref{theo_relative_QE_DCVF}.
	\end{proof}

	\begin{prop}\label{prop_geom_sorts_pur_im}
		In the theory $\mathrm{DCVF}$, the geometric sorts other than $K$ are blind to $K$.
		
	\end{prop}

	\begin{proof}
		We use additivity of transcendence degree, and the definable types in lattices and open balls that we have already encountered. Let $A$ be a set of field elements, $u$ a finite tuple in $\bigcup_n S_n \cup \bigcup_n T_n$, and $c$ an element of the field, such that $c \in acl(Au)$. We have to show that $c$ is in $acl(A)$. Without loss of generality, we may assume that there is an $n$ such that $u$ is of the form $([\Lambda_1], ..., [\Lambda_k], \alpha_1, ..., \alpha_k)$, where $[\Lambda_i] \in S_n$ and $\alpha_i \in \Lambda_i/ \mathfrak{m} \Lambda_i$, for all $i$. We may also assume that $A$ is a differential field of parameters.

		Let $d$ be a tuple of constants lifting $[\Lambda_1], ..., [\Lambda_k]$ generically over $Ac u$. More precisely, we can pick $d$ realizing the tensor product $(\eta^c_{B\Lambda_1} \otimes \cdots \otimes \eta^c_{B\Lambda_k})|_{Acu}$. By the remark above, the tensor products commute here. Then, the bases of the lattices $\Lambda_i$ given by $d$ induce $k$-bases of the quotients $\Lambda_i/ \mathfrak{m} \Lambda_i$, i.e. isomorphisms $\Lambda_i/ \mathfrak{m} \Lambda_i \simeq k^n$. Through these isomorphisms, the $\alpha_i$ can be identified with coordinates in $k$, i.e. translates of the open ball $\mathfrak{m}$. Then, let $e$ be a tuple of constants lifting the coordinates of the $\alpha_i$ generically over $Ac d u$. More precisely, realize one tensor product of the generics of the open balls considered (here, the tensor products do not commute, but it will not matter for our purposes), over the parameters $A c d u$.

		As all the generic types involved are transcendental, and $|e| = |d| = k n^2$, the transcendence degree of $d \flex e$ over $A$ is $2kn^2$. Moreover, since $c$ is in $acl(A de)$, and $d, e$ are made of constants, Corollary \ref{coro_closures_in_dcvf} above implies that $c$ is in the field-theoretic algebraic closure of $A (d,e)$. So, the transcendence degree of $c \flex d \flex e$ over $A$ is $2kn^2$. On the other hand, $c$ is not algebraic over $A$, and $d\flex e$ is algebraically independent over $Ac$, so the trancendence degree of $c \flex d \flex e$ over $A$ is $1+ 2k n^2$, which is a contradiction.\end{proof}

\end{subsubsection}

\begin{subsubsection}{Groups and fields interpretable in DCVF}

	\begin{prop}\label{prop_DCVF_morphism_with_imaginary_kernel}
		Let $G$ be an interpretable group. Assume that it is definably amenable. Then there exists an algebraic group $H$ and a definable morphism $G^{00} \rightarrow H$ whose kernel is purely imaginary.
	\end{prop}

	\begin{proof}
		
		We wish to apply Theorem \ref{theo_groups_in_algebraically_bounded_differential_fields}. By Proposition \ref{prop_geom_sorts_pur_im}, all geometric sorts different from the field sort are purely imaginary. Also, by Corollary \ref{coro_closures_in_dcvf}, the theory DCVF is a theory of algebraically bounded differential fields. Finally, by the fourth point of Fact \ref{fact_DCVF}, elimination of imaginaries holds in the geometric sorts, and the theory is NIP. \end{proof}

	\begin{prop}\label{prop_DCVF_def_subfields}
		Let $F$ be a definable subfield of $K \models DCVF$. Then, either $F$ is the field of constants of $K$, or $F$ is $K$ itself.
		
	\end{prop}

	\begin{proof}
		Let us first deal with the case where $F$ is of finite rank, i.e. there is a common differential equation satisfied by all the elements of $F$. Then, the Kolchin closure of $F$ is a $DCF$-definable subring of $K$ of finite rank, so it is the field of constants. Thus, $F$ is a subfield of the field of constants. Now, using the description of definable sets in $\mathrm{DCVF}$, we know that $F$ contains a nonempty valuative open subset of the constants. So there exists an open ball $b$ centered at $0$ such that $F$ contains all the constants in $b$. Then, if $x \in K^\partial$, there exists $a \in (b \cap K^\partial) \setminus \lbrace 0 \rbrace$ such that $a \cdot x \in b \cap K^\partial$. So $x = (a \cdot x) \cdot a^{-1} \in F$, and $F=K^\partial$, as required.

		Then, let us consider the case where $F$ is of infinite rank, i.e. contains a differentially transcendental element. Note that, by the previous case, the infinite definable field $F \cap K^\partial$ is the field of constants, i.e. $F$ contains the field of constants. On the other hand, by the description of definable sets in $\mathrm{DCVF}$, and the assumption of differential transcendence, there exists a nonempty valuative open subset $U$ of some $K^n$ such that $F$ contains the set $X := \lbrace x \in K \, | \, \nabla_n x \in U \rbrace$. Since $U$ is open, it contains a product of balls $b_0 \times \cdots \times b_{n-1}$. Up to translating $X$ by some element of $X$, we may assume that all these balls contain $0$. Now, let $a \in K$ be arbitrary. Since there exist arbitrarily small constants in $K$, there exists $\varepsilon \in K^\partial \setminus \lbrace 0 \rbrace$ such that $a \cdot \varepsilon \in X \subseteq F$. Thus, both $a \cdot \varepsilon$ and $\varepsilon$ are in $F$, so $a$ is also in $F$, and $F=K$.
	\end{proof}

	We now have enough tools to classify the interpretable fields.

	\begin{theo}\label{theo_fields_interpretable_in_DCVF}
		Let $F$ be an infinite field interpretable in $DCVF^{eq}$. Then, $F$ is definably isomorphic to either the residue field, the valued field, or the field of constants.
	\end{theo}
	
	\begin{proof}
		Assume that $F$ is not definably isomorphic to the residue field. Then, by relative quantifier elimination for DCVF, and the classification of interpretable fields in ACVF (Theorem \ref{theo_fields_ACVF}), we know that $F$ is not purely imaginary, i.e. $F$ is not blind to the field sort. Then, by Theorem \ref{theo_groups_in_algebraically_bounded_differential_fields}, the field $F$ admits a definable embedding into some finite extension of the valued field $K$. Since the latter is algebraically closed, this means that $F$ is definably isomorphic to a definable subfield of $K$. We can then conclude using Proposition \ref{prop_DCVF_def_subfields}.
	\end{proof}

\end{subsubsection}

\end{subsection}

\begin{rem}\label{rem_what_about_examples}
There are several other examples of algebraically bounded differential fields. In fact, all \emph{large geometric fields equipped with a generic derivation} are such. See \cite{ppp_def_groups_codf}[Fact 5.7] for the definition, and \cite{ppp_def_groups_codf}[Lemma 5.9] for the computation of the algebraic closure.
In such fields, if one can classify the  \emph{definable subfields} of the ambient field, one can deduce from Theorem \ref{theo_groups_in_algebraically_bounded_differential_fields} a classfication of the \emph{definable fields}. Then, if one can also classify the imaginaries, show that they are blind to the field sort, and classify the purely imaginary fields, one gets the classification of all interpretable fields.

For instance, the theory CODF of closed ordered differential fields, which is the model completion of the theory of ordered differential rings, eliminates imaginaries down to the field, and, as above, the only nontrivial definable subfield of the ambient field is the field of constants, which is also a real-closed field. Thus, the only interpretable fields are the field of constants, the field itself, and their algebraic closures.
\end{rem}

%
%
%
%
%
%
%
%
%
%
%
%
%
%

\begin{subsection}{Henselian valued fields}
	
	In this subsection, we let $Hen_0$ denote the theory of non-trivial henselian valued fields of characteristic $0$, in the three-sorted language $(K, k, \Gamma)$.

	\begin{defi}\label{defi_fin_rami}
		A valued field $(K, v)$ is \emph{finitely ramified} if it is of equicharacteristic $0$, or its residue characteristic is a prime number $p$, and the interval $[0, v(p)] \subseteq \Gamma$ is finite.
	\end{defi}

	
	\begin{prop}\label{prop_Hen_0}
		Let $T \supseteq Hen_0$ be a theory of finitely ramified henselian valued fields of characteristic $0$, in the three-sorted language $(K, k, \Gamma)$. Then, the following properties hold: 
		\begin{enumerate}
			\item The theory $T$ is an algebraically bounded theory of fields: the valued field sort part of the model-theoretic algebraic closure of any set $A$ of valued field elements is equal to the field-theoretic relative algebraic closure of the field generated by $A$, for any $A$ contained in a model $M \models T$.
			\item The residue field $k$ and the value group $\Gamma$ are stably embedded and orthogonal. Their induced structures are, up to naming constants, that of a pure field and a pure ordered abelian group respectively.
			\item The residue field and value group are blind to the valued field.
		\end{enumerate}

	\end{prop}
	
	\begin{proof}
		For the second point, see  \cite[Theorem 6.2]{Anscombe2023AxKochenErshovPF}. The first point is folklore, and can be proved using the elimination of field quantifiers in the RV language (see \cite[Theorem B]{Bas-EQHens}). Let us now prove the third point. Let $M \models T$ be sufficiently saturated. Let $a$ be an element of $K(M)$, $b \subset  K(M)$ and $\gamma \subset k(M) \cup \Gamma(M)$ be finite tuples. Assume that $a \in acl(b\gamma)$, and let us show that $a \in acl(b)$. 
		
		\begin{claim}
			There exists a finite tuple $c \in K(M)$ such that $\gamma \in dcl(c)$, and $c$ is an algebraically independent tuple, in the sense of field transcendence, over $a b$.
			
		\end{claim}
		\begin{proof}
			This can be proved by induction on the length of $\gamma$, using compactness. To lift one element of $k(M)$, recall that any non-trivial open ball is infinite, thus, by saturation, admits a point in $K(M)$ which is transcendental over $ab$. Then, add such a point to the tuple $b$, and proceed. To lift one element of $\Gamma(M)$, use the fact that any annulus contains a non-trivial open ball, and apply the same compactness argument as before.
		\end{proof}

		Then, using the claim, it suffices to compute transcendence degrees. By the first point of the proposition, since $a \in acl(b\gamma)\subseteq acl(bc)$, we have $trdeg(a/bc)=0$. By construction of $c$, we know that $trdeg(c/ab) = |\gamma| = trdeg(c/b)$. Then, by additivity, we have $trdeg(c / a b) + trdeg(a / b) = trdeg(c a / b) = trdeg(c/b) + trdeg(a/c b) = trdeg(c/b)$. Thus, we have $trdeg(c / a b) + trdeg(a / b) = trdeg(c / b)$. Since $trdeg(c / ab) = trdeg(c/b)$, this implies $trdeg(a/b) = 0$, as required.
	\end{proof}

	For equicharacteristic zero henselian valued fields, recent results (in \cite{RK-Vic-Imaginaries}) give a precise description of imaginary sorts; the third point of Proposition \ref{prop_Hen_0} holds for all sorts different from the valued field sort, in a language which admits elimination of imaginaries. Before we can state the results, we need to introduce some definitions.
	
	\medskip \medskip
	
	For the remainder of this subsection, we fix a sufficiently saturated equicharacteristic $0$ henselian valued field $(K, v)$, and let $M^{eq}$ denote the underlying structure in the imaginary sorts. We let $\mathrm{Cut}$ denote the $\vee$-definable set of codes of cuts in $\Gamma$, i.e. upwards closed definable subsets of $\Gamma$. We also let $\mathrm{Cut}^*$ denote $\mathrm{Cut} \setminus \lbrace \emptyset, \Gamma \rbrace$.
	We assume that $\Gamma$ is of \emph{bounded regular rank}, i.e. has at most countably many definable convex subgroups.
	
	If $\Delta \leq \Gamma$ is a definable convex subgroup, we let $v_{\Delta} : \Gamma \rightarrow \Gamma / \Delta$ denote the definable quotient map.
	
	For any $c \in \mathrm{Cut}(M)$, and any $a \in K(M)$, we let $b_c(a)$ denote the \emph{generalized ball} of cut $c$ around $a$, defined by $b_c(a) = \lbrace x \in K \, | \, v(x-a) \in c \rbrace$. We let $\mathrm{B}_g$ denote the $\vee$-definable set of codes of generalized balls. In particular, if $a=0$, we define $I_c = b_c(0)$; it is a sub-$\mathcal{O}$-module of $K$. 
	
	For each finite tuple $c=(c_i)_i$ in $\mathrm{Cut}^*$, we define a set of (codes of) $\mathcal{O}$-modules of type $c$, which we denote $\mathrm{Mod}_c$. More precisely, let $B_n$ denote the group of invertible upper triangular matrices, and consider the map $\mu_c$ which maps a matrix, i.e. a triangular basis $(a_i)_{i <n}$ of $K^n$, to the module $\sum_i I_{c_i} a_i$. Then, let $\mathrm{Mod}_c$ denote its image, so that we have a surjective definable map $\mu_c  : B_n \rightarrow \mathrm{Mod}_c$. Note that the $\mathrm{Mod}_c$ vary definably in $c$, so, for each $n$, we get a $\vee$-definable set $\mathrm{Mod}_n$. We let $\mathrm{Mod}$ denote the $\vee$-definable set $\bigsqcup_n \mathrm{Mod}_n$. The restriction to $\mathrm{Cut}^*$ is justified by the following:
	
	\begin{fact}{(see Corollary 3.7 in \cite{RK-Vic-Imaginaries})}
		Any $M$-definable sub-$\mathcal{O}$-module of $K^n$ can be encoded in $K \cup \mathrm{Mod}^{eq}$.
	\end{fact}

	The following is a weak consequence of Proposition 6.1 in \cite{RK-Vic-Imaginaries}:
	\begin{prop}\label{prop_weak_elim_hen_00}
		Let $e \in M^{eq}$ and $A = acl(e)$. Then $e \in dcl((K \cup \mathrm{Mod}^{eq})(A))$.
	\end{prop}
	\begin{proof}
		The result  \cite[Proposition 6.1]{RK-Vic-Imaginaries} yields that $e$ is definable over $\mathcal{G}(A) \cup (\mathrm{RV} \cup \mathrm{Lin}_A)^{eq} (A)$, where $\mathcal{G}$ denotes the geometric sorts, and $\mathrm{Lin}_A$ the $A$-definable ``$k$-linear sorts''. The point is that all the sorts in $\mathcal{G}$, RV and $\mathrm{Lin}_A$ can be embedded $A$-definably in products of sorts in $K \cup \mathrm{Mod}^{eq}$. \end{proof}

	The other result of interest to us is the following unary decomposition:

	\begin{prop}{(see Proposition 3.9 in \cite{RK-Vic-Imaginaries})}\label{prop_unary_dec}
		Let $s \in \mathrm{Mod}_c$, where $c$ is a tuple in $\mathrm{Cut}^*$, with $|c| = n$. There exists a tuple $b= (b_l)_l \in M^{eq}$ (each $b_l$ being identified with a subset of some $K^{r_l}$) and a family of interpretable sets $(X_l)_l$, where each $X_l$ is definable over $cb_{<l}$, such that: 
		
		\begin{itemize}[$\bullet$]
			\item For every $l$, we have $b_l \in X_l$.
			
			\item We have $dcl(s) = dcl(cb)$.
			
			\item If $l < n$, we have $X_l = \Gamma/\Delta_{c_i}$, where $l= (i,i)$.
			
			\item If $l \geq n$, then $X_l$ has a $cb_{<l}$-definable $K/I_l$-torsor structure, where $l=(i,j)$, and $I_l$ is a $cb_{<l}$-definable multiple of the module $(I_{c_i} : I_{c_j}) = \lbrace x \, | \, x I_{c_j} \subseteq I_{c_i} \rbrace$.
			
			Moreover, for any choice of $(a_l)_l$ with $a_l \in b_l$ for all $l$, there is a family of (uniformly) $ca_{<l}$-definable isomorphisms $f_l : X_l \rightarrow K/I_l$ and a $ca_{<l}$-definable function $g_l : f_l(b_l) \rightarrow b_l$.
			
		\end{itemize}
	\end{prop}
	
	
	%
	%
	%
	%
	%
	
	From that, we deduce that the sorts in $\mathrm{Mod}$ are blind to the field sort:

	\begin{prop}\label{prop_Mod_blind_to_K}
		For all $m$, if $Y \subseteq \mathrm{Mod}_m$ is a definable set, then it is blind to $K$.
	\end{prop}
	
	\begin{proof}
		By compactness, it suffices to prove that, for all $m$, for all $(c_i)_{i < m}$ in $\mathrm{Cut}^*$, for all $s \in \mathrm{Mod}_c$, for all parameter sets $A_0 \subseteq K(M)$, for all $a \in K$, if $a \in acl(A_0 cs)$, then $a \in acl(A_0)$. Let us now prove it. 
		
		
		We now use Proposition \ref{prop_unary_dec}. Let $(b_l)_l$ and $(X_l)_l$ be as in the conclusion of Proposition \ref{prop_unary_dec}. 
		
		\begin{claim}
			There exists a finite tuple of field elements $d$ such that $d$ is algebraically independent from $a$ over $A_{0, K}$,  the field generated by $A_0$, and $cb \in dcl(d)$.
		\end{claim}
		
		\begin{proof}
			First, recall that any element of $\Gamma$ has infinitely many preimages in $K$, and that $\Gamma$ is stably embedded. So, by compactness, we may find appropriate lifts in $K$, i.e. preimages under $0$-definable maps,  of both the $c_i$ and the $b_l$, for $l < n$. Let $d_{<n}$ denote such a finite tuple in $K$, i.e. $d_{<n}$ is algebraically independent over $A_{0, K} a$, and $cb_{<n} \in dcl(d_{<n})$. Now, it remains to deal with the $K/I_l$-torsors $X_l$, for $l \geq n$. 
			Proceeding by induction, assume we have a suitable tuple $d_{<l}$ lifting $cb_{<l}$. Then, since $d_{<l}$ is algebraically independent from $a$ over $A_{0, K}$, using algebraic boundedness, we can find a model $M_l$ containing $d_{<l}$, such that $M_l$ is algebraically independent from $a$ over  $A_{0, K}$: following the proof of the Löwenheim-Skolem theorem, using the Tarski-Vaught criterion, if a formula has finitely many solutions, any of them is algebraic over the current set; if it has infinitely many solutions, by compactness, it has one which is transcendental over all the parameters. Then, over $M_l$, there is an isomorphism $X_l \simeq K/I_l$. Let $a_l \in M_l$ be a finite tuple over which this isomorphism is defined. Since each coset in $K/I_l$ has infinitely many preimages in $K$, for $I_l$ is nontrivial, we can, by compactness, find a lift $e_l \in K$ of (the image of) $b_l$ such that $e_l$ is transcendental over $A_0 a M_l$. Then, put $d_l = a_l e_l$, and $d_{\leq l} = d_{<l} d_l$. This concludes the proof.
		\end{proof}
		Then, recall that $dcl(s) = dcl(cb)$. So, we have $a \in acl(A_0 cs) \subseteq acl(A_0 d)$. Thus, since pure henselian valued fields of characteristic $0$ are algebraically bounded, the element $a$ is algebraic, in the field sense, over $A_{0, K}(d)$. Since $d$ is algebraically independent from $a$ over $A_{0, K}$, this implies that $a$ is algebraic over $A_0$, which concludes the proof.
	\end{proof}
	
	We can now prove the following
	
	\begin{prop}\label{prop_good_coords_hen_00}
		Let $X$ be a definable set in $M^{eq}$. Then, there exists a definable set $Y$, an integer $n$ and a definable finite-to-one correspondence $X \rightarrow Y \times K^n$, such that $Y$ is blind to $K$.

	\end{prop}

	\begin{proof}
		By Proposition \ref{prop_weak_elim_hen_00}, we have weak elimination of imaginaries down to $K \cup \mathrm{Mod}^{eq}$. Since we are working up to finite correspondence, and since blindness to $K$ is preserved by taking quotients, it suffices to prove the result in the case where $X~\subseteq~K^n \times~{(\mathrm{Mod}_m)}^d$, for some $n, m, d$. 
		The result then follows from Proposition \ref{prop_Mod_blind_to_K}.
	\end{proof}

	\begin{coro}\label{coro_good_coords_hen_00} Let $M$ be a model of a theory of pure henselian valued fields of equicharacteristic $0$ whose value group is of bounded regular rank. 
		Let $X$ be a definable set in $M^{eq}$. Then, there exists an integer $N$ and a definable map $X \rightarrow K^N$, whose fibers are blind to $K$.

	\end{coro}
	\begin{proof}
		Using Proposition \ref{prop_good_coords_hen_00}, we get a definable finite-to-one correspondence  $X \rightarrow Y \times K^n$, where $Y$ is blind to $K$. Then, composing with projection to $K^n$, and using elimination of finite sets in fields, we get a definable map $X \rightarrow K^N$, whose fibers admit a definable finite-to-one correspondence with some finite disjoint union of copies of $Y^n$. Since $Y$ is blind to $K$, these fibers are as well.
	\end{proof}

	\begin{theo}\label{theo_groups_hen_00}
		Let $T$ be a theory of pure henselian valued fields of equicharacteristic $0$, whose value groups have bounded regular rank, i.e. at most countably many definable convex subgroups. Assume that $T$ is NIP, or that $T$ is $\mathrm{NTP_2}$ with $Gal(K^{alg}/K)$ bounded in all models of $T$. Let $M$ be a sufficiently saturated model of $T$, and $G$ be an $M$-definable definably amenable group. Then, there exists an algebraic group $H$ over $K(M)$ and an $M$-definable group homomorphism $G^{00}_M \rightarrow H$ whose kernel is blind to the field sort.
		
	\end{theo}
	
	\begin{proof}
		First, by Corollary \ref{coro_good_coords_hen_00}, there exists an integer $N$ and an $M$-definable map $G \rightarrow K^N$, whose fibers are blind to the field sort $K$. Also recall that, by Proposition \ref{prop_Hen_0}, the theory $T$ is algebraically bounded: model-theoretic and field-theoretic algebraic closure coincide in the valued field sort. Then, it suffices to apply Theorem \ref{theo_general_morphism_with_imaginary_kernel} to conclude.
	\end{proof}

	\begin{theo}\label{theo_fields_hen_0}
		Let $T \supseteq Hen_0$ be a theory of finitely ramified pure henselian valued fields of characteristic $0$ with perfect residue fields, in the three-sorted language. Assume that $T$ is NIP, or that $T$ is $\mathrm{NTP_2}$ with $Gal(K^{alg}/K)$ bounded in all models of $T$. Also assume that $k$, with its induced structure, is algebraically bounded as a field. Let $M$ be a sufficiently saturated model of $T$. Let $m, n, l < \omega$. Let $F \subseteq K^n \times k^m \times \Gamma^l$ be an $M$-definable infinite field. Then, exactly one of the following holds:
		
		\begin{enumerate}
			
			\item The field $F$ admits a definable embedding into some finite extension of $K$.
			
			\item The field $F$ admits a definable embedding into some finite extension of $k$.
			
			\item The field $F$ is definably isomorphic to a field definable in $\Gamma$.
			
		\end{enumerate}

	\end{theo}

	\begin{proof}
		By Proposition \ref{prop_Hen_0}, the theory $T$ is algebraically bounded, the sorts $k, \Gamma$ are orthogonal and blind to $K$, and stably embedded.
		Note that the conclusions are pairwise inconsistent, because $k$ and $\Gamma$ are orthogonal, and both are blind to $K$.

		Assume that $F$ does not admit a definable embedding into a finite extension of $K$. Then, by Corollary \ref{coro_fields_are_either_pur_im_or_algebraic}, the field $F$ is purely imaginary, which means that its projection to $K^n$ is finite. Thus, encoding a finite set, and changing $m, l$, if necessary, we may assume that $F \subseteq k^m \times \Gamma^l$. Note that, by stable embeddedness (see \cite[Theorem 6.2]{Anscombe2023AxKochenErshovPF}), the field $F$ is definable in the induced structure $(k(M), \Gamma(M))$, which is sufficiently saturated. Also note that, in the $\mathrm{NTP_2}$ case, by henselianity and perfection of $k$, the Galois group $Gal(k^{alg}/k)$ is also bounded. Then, the hypotheses imply that Corollary \ref{coro_fields_are_either_pur_im_or_algebraic} can be applied to the theory of the pair $(k, \Gamma)$. So, either $F$ admits a definable embedding into  some finite extension of $k$, or $F$ can be embedded, as a definable set, into some power of $\Gamma$. Since $\Gamma$ is stably embedded, the latter case implies that $F$ is definably isomorphic to a field definable in $\Gamma$. This concludes the proof.
	\end{proof}
	%
	%
	
	\begin{rem}
		\begin{enumerate}
			\item It seems likely that there are no infinite definable fields in $\Gamma$, for instance because of the result \cite[Corollary 1.10]{CluHalQEOAG}, which states that, in pure ordered abelian groups, definable maps are piecewise linear. However, proving such results is beyond the scope of this paper.
			
			\item The hypothesis of finite ramification is not essential. The main ingredients are orthogonality between $k$ and $\Gamma$, and blindness of these sorts to $K$. 
		\end{enumerate}
	\end{rem}

	\begin{ex}
		The theorem applies to any {pure} finitely ramified henselian valued field $(K, v)$ of characteristic $0$, whose residue field is elementary equivalent to either $\mathbb{C}$, $\mathbb{R}$, $\mathbb{F}_{p}^{alg}$, a finite extension of $\mathbb{Q}_p$, or a finite field. 
	\end{ex}

\end{subsection}

\end{section}

\printbibliography[
heading=bibintoc,
title={References}
]

\end{document}